\documentclass{amsart}
\usepackage{amsmath,amsthm}
\usepackage{graphics}
\usepackage{color}
\usepackage{epsfig}

\graphicspath{{./Figures/}}

\theoremstyle{plain}
\newtheorem{thm}{Theorem}[section]
\newtheorem{cor}[thm]{Corollary}

\newtheorem{lemma}[thm]{Lemma}
\newtheorem{claim}[thm]{Claim}

\newtheorem{question}[thm]{Question}

\theoremstyle{definition}
\newtheorem{remark}[thm]{Remark}

\newcommand{\comment}[1]{}

\newcommand{\bdry}{\ensuremath{\partial}}
\DeclareMathOperator{\Int}{Int}

\DeclareMathOperator{\lcm}{lcm}

\newcommand{\R}{\ensuremath{\mathbb{R}}}
\newcommand{\Z}{\ensuremath{\mathbb{Z}}}

\newcommand{\opt}{{\small OPT}}
\newcommand{\gof}{{\small GOF}}

\begin{document}

\title{Once-punctured tori and knots in lens spaces.}
\author{Kenneth L.\ Baker} 
\address{School of Mathematics, Georgia Institute of Technology \\ Atlanta, Georgia 30332}
\email{kb@math.gatech.edu}

\begin{abstract}
We determine the non-null homologous knots in lens spaces whose exteriors contain properly embedded once-punctured tori.  All such knots arise as surgeries on the Whitehead link and are grid number $1$ in their lens spaces.  As a corollary, we classify once-punctured torus bundles that admit a lens space filling.
\end{abstract}

\maketitle

\section{Introduction}

Say a knot whose exterior contains a properly embedded once-punctured torus is an {\em \opt--knot}.  In this article we determine every non-null homologous \opt--knot $K$ in a lens space $X$ up to homeomorphism of the pair $(X,K)$.  As a corollary, when combined with previous work of the author \cite{baker:cgofkils}, we determine all once-punctured torus bundles that admit a lens space filling.  See Theorem~\ref{thm:main} and Corollary~\ref{cor:optbundle} respectively.   

The lens space $L(r,q)$ is obtained by $-r/q$ Dehn surgery on the unknot in $S^3$ for coprime integers $r$ and $q$.  Recall that $L(r,q) \cong L(r',q')$ if and only if $r'=\pm r$ and $q' = \pm q^{\pm 1}$ modulo $r$.

Let $K$ be an \opt--knot in a lens space $X$ and let $S$ be a once-punctured torus properly embedded in its exterior $E(K) = X - N(K)$.  Let $s$ be the minimal positive integer for which $s[K] = 0$ in $H_1(X)$.  We say $s$ is the {\em order} of $K$.  Observe that $s$ is the number of times a meridian of $K$ minimally intersects $\bdry S$ on $\bdry E(K)$.  Let $r$ be the order of $X \not \cong S^1 \times S^2$ so that $H_1(X) \cong \Z/r$, $X \cong L(r,q)$ for some $q$, and $s \leq r$.  If $X \cong S^1 \times S^2$, then $H_1(X) = \Z$, we take $r = 0$, and allow $s=0$.

\begin{figure}[h]
\centering
\input{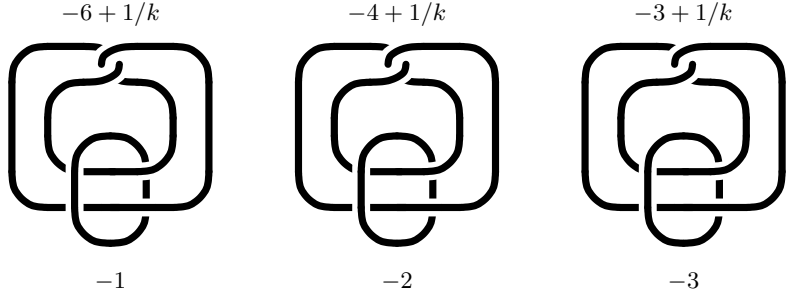}
\caption{Surgery descriptions of the five knotted families of non-null homologous \opt--knots in lens spaces.} 
\label{fig:whiteheadsurgdescr}
\end{figure}

Let $W$ denote the Whitehead link depicted three times in Figure~\ref{fig:whiteheadsurgdescr}.  Let $W(\gamma_1, \gamma_2)$ denote Dehn surgeries of slopes $\gamma_1$ and $\gamma_2$ on the two components of $W$ in the standard coordinates.  We say a knot $K$ is {\em unknotted} if its exterior is a solid torus 
and is {\em trivial} if it is the boundary of an embedded disk. 
The {\em type} of a torus knot is defined in \S\ref{sec:torusknots}.

\begin{thm}\label{thm:main}
If $K$ is a non-null homologous \opt--knot in a lens space then up to homeomorphism for some integer $k \neq 0$ either
\begin{itemize}
\item[(I)] $K$ is the type $\{2,3\}$ torus knot in $L(6k-1,2k-1)$ arising as the core of the $-6+1/k$ surgery in $W(-1,-6+1/k)$, 
\item[(II)] $K$ is the type $\{2,4\}$ torus knot in $L(8k-2,4k+1)$ arising as the core of the $-4+1/k$ surgery in $W(-2,-4+1/k)$, 
\item[(III)] $K$ is the type $\{3,3\}$ torus knot in $L(9k-3,3k-2)$ arising as the core of the $-3+1/k$ surgery in $W(-3,-3+1/k)$, 
\item[(IV)] $K$ is the knot in $L(9k-3,3k-2)$ arising as the core of the $-3$ surgery in $W(-3,-3+1/k)$,
\item[(V)] $K$ is the knot in $L(8k-2,4k+1)$ arising as the core of the $-2$ surgery $W(-2,-4+1/k)$, or
\item[(VI)] $K$ is an unknotted knot in $L(r,q)$ arising as the core of the $-r/r$ surgery in $W(-r/q, \infty)$, $|r|\neq1$.
\end{itemize}
See Figure~\ref{fig:whiteheadsurgdescr}.  The once-punctured torus compresses in the exterior of $K$ if and only if $K$ belongs to family (VI).
\end{thm}

\begin{proof}
If the exterior of $K$ is a solid torus, then it is the core of $-r/q$ Dehn surgery on the unknot in $S^3$ and belongs to family (VI).  Now assume the exterior of $K$ is not a solid torus. When $K$ is a torus knot then by Lemma~\ref{lem:torusknots} it belongs to one of the first three families.  When $K$ is not (necesarily) a torus knot then for $s \geq 3$ it follows from Theorem~\ref{thm:sgeq3} that $K$ belongs to family (IV).  For $s=2$ it follows from Theorem~\ref{thm:seq2} that $K$ belongs to family (V).  
\end{proof}

\begin{remark}
Note that the first three families of knots are all torus knots of orders $|-6k+1|$, $|-4k+1|$, and $|-3k+1|$ respectively.  The fourth and fifth families of knots are not torus knots (except when $k=1$) and their orders are $3$ and $2$ respectively.  Observe that, consequentially, only the knots in the first family represent generators of homology in their lens spaces.
\end{remark}

A (knotted) knot or link in a $3$--manifold has {\em tunnel number $1$} if the complement of an open neighborhood of some arc properly embedded in the exterior of the knot or link is a handlebody.  Unknotted knots have tunnel number $0$.

\begin{cor}
Each knotted non-null homologous \opt--knot in a lens space has tunnel number $1$.
\end{cor}

\begin{proof}
Theorem~\ref{thm:main} shows that each non-null homologous \opt--knot in a lens space arises from Dehn surgery on the Whitehead link which itself has tunnel number $1$.
\end{proof}

\begin{cor}
For each $k \in \Z$, $k \neq 0$, the lens space $L(6k-1,2k-1)$ contains exactly one knotted non-null homologous \opt--knot and the lens spaces $L(8k-2,4k+1)$ and $L(9k-3,3k-2)$ each contain exactly two knotted non-null homologous \opt--knots.
\end{cor}

\begin{proof}
From Theorem~\ref{thm:main} it is clear that these lens spaces contain at least as many as claimed.  To show they contain no more, we must examine when these three families of lens spaces coincide.  First observe that $2 \vert 8k-2$ and $3 \vert 9k-3$ but $2 \!\!\! \not \vert 6k-1$ and $3 \!\!\!\not \vert 6k-1$; hence $L(6k-1,2k-1)$ contains no other non-null homologous \opt--knots.  Now notice that since $(3k-2)^2 \equiv 1$ modulo $9k-3$, if $L(9k-3,3k-2) \cong L(8\ell-2,4\ell+1)$ for some $k, \ell \in \Z$ then, from the classification of lens spaces, it must be the case that $(4\ell+1)^2 \equiv 1$ modulo $8\ell-2$.  But this only occurs if $\ell = 1$, and hence $k=1$ too.  In this case the surgery descriptions of the knots in $L(6,1)\cong L(6,5)$ coincide as $W(-2,-3) = W(-3,-2)$
\end{proof}

\subsection{Lens space fillings of once-punctured torus bundles.}
By a {\em genus one fibered knot}, or {\em \gof--knot} for short, we mean a null homologous knot in a closed orientable $3$--manifold whose exterior is a once-punctured torus bundle; a \gof--knot is the binding of an open book whose pages are once-punctured tori.

Burde and Zieschang use the homology of once-punctured torus bundles in their proof that a \gof--knot in $S^3$ is either a trefoil or has exterior homeomorphic to that of the figure $8$, \cite{bz:nkuf}.  Gonz\'ales-Acu\~na pins this down concluding that the \gof--knots in $S^3$ are the trefoils and the figure $8$, \cite{ga:dcok}.  

Passing from $S^3$ to lens spaces, Morimoto describes all the \gof--knots in a handful of lens spaces, \cite{morimoto:gofkils}.  In \cite{baker:cgofkils}, we extend Morimoto's result by classifying all the \gof--knots in all lens spaces.  Using the homology of once-punctured torus bundles in a manner similar to Burde and Zieschang, Baldwin observes that if a once-punctured torus bundle is the exterior of a non-null homologous knot in a lens space of prime order, then the once-punctured torus bundle is homeomorphic to the exterior of a trefoil, \cite{baldwin:anogofkils}.  We now complete the classification of knots in lens spaces with once-punctured torus bundle exteriors in Theorem~\ref{thm:fibered} and describe the monodromy of a once-punctured torus bundle admitting a lens space filling in Corollary~\ref{cor:optbundle}.

Let $x$ and $y$ be two simple closed curves on $S$ that transversally intersect once; let $\tau_x$ and $\tau_y$ be right-handed Dehn twists along these curves.

\begin{thm}\label{thm:fibered}
If $K$ is a non-null homologous knot in the lens space $X$ with a once-punctured torus bundle of monodromy $\phi$ as its exterior, then up to homeomorphism for some integer $k\neq0$ either 
\begin{itemize}
\item[(I)] $K$ is a type $\{2,3\}$ torus knot in $L(6k-1,2k-1)$ with $\phi \cong \tau_x \tau_y$,
\item[(II)] $K$ is a type $\{2,4\}$ torus knot in $L(8k-2,4k+1)$ with $\phi \cong \tau_x^2 \tau_y$, 
\item[(III)] $K$ is a type $\{3,3\}$ torus knot in $L(9k-3,3k-2)$ with $\phi \cong \tau_x^3 \tau_y$, 
\item[$\underset{k=-1}{\mbox{(IV)}}$]  $K$ is the non-torus knot in $L(12,5)$ arising as the core of the $-3$ surgery in $W(-3,-4)$ with $\phi \cong \tau_x^4 \tau_y$, or
\item [$\underset{k=-1}{\mbox{(V)}}$] $K$ is the non-torus knot in $L(10,3)$ arising as the core of the $-2$ surgery in $W(-2,-5)$ with $\phi \cong \tau_x^5 \tau_y$. \end{itemize}
\end{thm}
\begin{proof}
When $K$ is a torus knot families (I), (II), and (III) follow from Lemma~\ref{lem:torusknots}.  When $s\geq3$ Theorem~\ref{thm:sgeq3} and in particular Lemma~\ref{lem:order3kplus1} gives the knot in family (IV) with $k=-1$; the knot in family (IV) with $k=+1$ is a torus knot.  When $s=2$ Theorem~\ref{thm:seq2} and in particular Lemma~\ref{lem:order2kplus1} gives the knot in family (V) with $k=-1$; the knot in family (V) with $k=+1$ is a torus knot.
\end{proof}

\begin{remark}
In Theorem~\ref{thm:fibered} the first three monodromies are periodic, the fourth is reducible, and the fifth is pseudo-Anosov (see e.g.\ \cite{cassonbleiler}).  Indeed the fourth exterior, $W(-4,\cdot)$, contains an essential Klein bottle.  The fifth exterior, $W(-5,\cdot)$, is the exterior of the figure $8$ sister.
\end{remark}

\begin{cor}\label{cor:optbundle}
Up to homeomorphism, a once-punctured torus bundle admits a lens space filling if and only if it has monodromy $\phi \cong \tau_x^k \tau_y^2 \tau_x^\ell \tau_y^{-1}$ for $k, \ell \in \Z$.
\end{cor}
\begin{proof}
Let $\Delta$ denote the distance of the slope of a lens space filling of a once-punctured torus bundle from the slope of the boundary of a fiber, see \S\ref{sec:slopes}.
 
If $\Delta >1$, then the core of the filling is a knotted non-null homologous \opt--knot  in a lens space with fibered exterior.  Hence $\phi$ must be as in Theorem~\ref{thm:fibered}.  These agree with choosing $k \in \{1,2,3,4,5\}$ and $\ell =0$ in the expression $\phi = \tau_x^k \tau_y^2 \tau_x^\ell \tau_y^{-1}$.

If $\Delta = 1$, then the corresponding knot is a (null homologous) genus one fibered knot.  In \cite{baker:cgofkils} we show that $\phi$ has the stated form and that every such $\phi$ is the monodromy of a genus one fibered knot in a lens space.

If $\Delta = 0$ then the filling is a torus bundle.  But this cannot be a lens space.
\end{proof}

\begin{cor}
Every non-null homologous knot $K$ in a lens space with a once-punctured torus bundle exterior admits a non-trivial lens space surgery.  In particular there is a \gof--knot in a lens space with the same exterior as $K$.
\end{cor}

\begin{proof}
The exteriors of the five families of non-null homologous knots in lens spaces  with once-punctured torus bundle exteriors listed in Theorem~\ref{thm:fibered} are $W(-1, \cdot)$, $W(-2,\cdot)$, $W(-3,\cdot)$, $W(-4, \cdot)$, and $W(-5, \cdot)$ respectively.  The core of the $\infty$ filling of the remaining component of $W$ gives a \gof--knot in $L(1,1)$, $L(2,1)$, $L(3,1)$, $L(4,1)$, and $L(5,1)$ respectively.
\end{proof}

\begin{remark}
Each of these \gof--knots in lens spaces that admit non-trivial lens space surgeries (begetting the knots of Theorem~\ref{thm:fibered}) may be visualized up to homeomorphism as the plumbing of the $+1$--Hopf band in $S^3$ onto the $r$--Hopf band in $L(r,1)$ for $r=1,2,3,4,5$.  We say an {\em $r$--Hopf link} is the two component fibered link in the lens space $L(r,1)$ with an annulus fiber whose monodromy is $r$ Dehn twists along the core curve of the annulus; i.e.\ the binding of an open book for $L(r,1)$ with annular pages.   An {\em $r$--Hopf band} is then an annulus fiber of an $r$--Hopf link.  See \cite{bjk:tnogofk}.
\end{remark}

\subsection{Grid number $1$ knots.}
Recently the advent of combinatorial link Floer homology for links in $S^3$, \cite{mos:acdokfh} and \cite{most:oclfh}, has brought to attention Cromwell's arc presentations of links,  \cite{cromwell:apokal} and \cite{dynnikov:apolms}.  These may be reinterpreted as toroidal grid diagrams of links.  One observes that analogous toroidal grid diagrams may be defined for knots in lens spaces too.  Grigsby employs this in \cite{grisgsby:cdokfhocbc}.   

We define {\em grid number $1$ knots} in section \S\ref{sec:gridnumber} and discuss them here.  Let us first note that in an oriented lens space the non-trivial elements of homology are canonically identified with the oriented grid number $1$ knots.  In various senses, these grid number $1$ knots are ``small'' or ``fundamental''.   For example, from the presentation of a knot $K$ as a grid number $1$ knot one quickly obtains that the Heegaard Floer knot homology of $K$ is exceptionally simple with just one generator in each Spin$^c$ structure.  In fact, each grid number $1$ knot is a U--knot, see Definition 5.2 of \cite{os:kfhars}.  Also, the known knots in $S^3$ that admit integral Dehn surgeries yielding lens spaces, i.e.\ the classical Berge knots, all correspond to grid number $1$ knots in the resulting lens spaces, see \cite{berge:skwsyls}.  Yet another aspect of the simplicity of grid number $1$ knots is demonstrated in the following theorem.

\begin{thm}\label{thm:gridnumber1}
Every non-null homologous \opt--knot in a lens space has grid number $1$.
\end{thm}
\begin{proof}
Refer to Theorem~\ref{thm:main} for the classification of non-null homologous \opt--knots in lens spaces.
That the grid number of $K$ is $1$ follows from Lemma~\ref{lem:torusknotgridnumber} for the torus knots,  from Lemma~\ref{lem:order3gridnumber} when $s \geq 3$, from Lemma~\ref{lem:order2gridnumber} when $s=2$, and trivially for the unknotted knots
\end{proof}

This suggests a question.  First, let us say a properly embedded, connected, orientable surface $S$ in the exterior $E(K) = Y - N(K)$ of a knot $K$ in some closed orientable $3$--manifold $Y$ is a {\em rational Seifert surface} for $K$ if an orientation on $S$ causes $\bdry S$ to be a coherently oriented collection of circles on $\bdry E(K)$.  Then the {\em (rational) genus of $K$} is the minimal genus among all rational Seifert surfaces for $K$. 
\begin{question}
Do the grid number $1$ knots in a lens space have minimal genus among knots in their homology class?
\end{question}
Note that since the number of boundary components of a rational Seifert surface is invariant for knots within a homology class, this question may be rephrased in terms of Euler characteristics.

Further supporting the simplicity of grid number $1$ knots in lens spaces, Ozsv\'ath and Szab\'o algebraically prove that a grid number $1$ knot in a lens space representing a generator of homology has fibered exterior, \cite{os:kfh+lss}.   They show this using a theorem of Stallings that gives an algebraic characterization of fibered orientable $3$--manifolds $Y$, \cite{stallings:ofc3m}, together with a theorem of Brown about the finite generation of kernels of characters of two-generator one-relator groups, \cite{brown:tvatbnsi}.  In accordance with this notion of simplicity, one might suspect that all grid number $1$ knots have fibered exterior.  Alas, this is not so.

\begin{cor}
Not every grid number $1$ knot has fibered exterior.
\end{cor}
\begin{proof}
By Theorem~\ref{thm:fibered}, families (IV) and (V) of non-null homologous \opt--knots are not fibered if $k \neq \pm1$.  This with Theorem~\ref{thm:gridnumber1} gives the result.
\end{proof}

Because the \opt--knots in families (II), (III), (IV), and (V) do not represent generators of homology (their orders are strictly less than the orders of the lens spaces that contain them) the application of Brown's theorem in concert with Stallings' theorem does not work as cleanly.  It is not uncommon, as in families (II) and (III) and in $L(p,1)$, to find grid number $1$ knots that are torus knots, and hence have fibered exterior, but do not represent generators of homology.  Together Theorem~\ref{thm:fibered} and Theorem~\ref{thm:gridnumber1} show that there are but two non-torus \opt--knots of grid number $1$ whose exteriors are fibered and yet do not represent generators of homology.

\begin{question}
Aside from torus knots and ones that represent generators of homology,
which grid number $1$ knots have fibered exterior?
\end{question}

\subsection{Acknowledgments}
This work was partially supported by NSF Grant DMS--0089927 and NSF Grant DMS--0239600.  We thank John Etnyre for some useful comments and suggestions.

\section{Preliminaries}

\begin{lemma}\label{lem:compress}
Let $S$ be a once-punctured torus properly embedded in the exterior of a non-null homologous knot $K$ in a lens space $X$.  If $S$ compresses then $K$ is an unknotted knot.
\end{lemma}
\begin{proof}
Let $S'$ be the disk resulting from a compression of $S$.  Then $N(K \cup S')$ is a punctured lens space of the same order as $K$.  Since lens spaces are irreducible, $X-N(K\cup S')$ is a $3$--ball.  Then $E(K) = X-N(K)$ is obtained by attaching a $1$--handle that is a collar neighborhood of $S'$ to this ball.  Hence $E(K)$ is a solid torus.
\end{proof}

\subsection{Slopes and surgeries}\label{sec:slopes}
A {\em slope} is an isotopy class of essential simple closed curves on a torus.
If $\gamma$ and $\delta$ are two slopes, then $\Delta(\gamma,\delta)$ denotes their {\em distance}, the minimal geometric intersection number among representatives of $\gamma$ and $\delta$ in the torus.  If $c$ and $d$ are elements of $\gamma$ and $\delta$ respectively then we also write $\Delta(c,d) = \Delta(\gamma, \delta)$.  For example, in this article we assume that for a meridian $\mu$ of $K$, $\Delta(\mu, \bdry S)=s$ on $\bdry N(K)$.  

Given a link $K = K_1 \cup \dots \cup K_n$, let $K(\gamma_1, \dots, \gamma_n)$ denote the manifold obtained by $\gamma_i$--Dehn surgery on each component $K_i$ for $i=1, \dots, n$.  We use the ``slope'' $\gamma_i = \cdot$ to indicate that a tubular neighborhood of $K_i$ is removed and left unfilled.

Let $W$ denote the Whitehead link shown in Figure~\ref{fig:whiteheadsurgdescr} and $N$ denote the pretzel link $P(-2,-2,-2)$.  The pretzel link $N=P(-2,-2,-2)$ has exterior homeomorphic to the the so-called ``magic'' manifold, see Gordon and Wu \cite{gordon:ssadf}, \cite{gordonwu:taadf} and, for our purposes in particular, Martelli and Petronio  \cite{mp:dfotm3m}.

\begin{figure}
\centering
\input{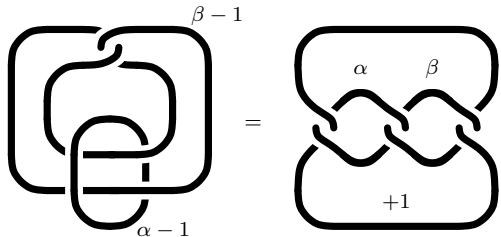}
\caption{Blowing down the $+1$ surgery gives the relationship $W(\alpha-1,\beta-1)=N(\alpha,\beta,+1)$. }
\label{fig:surgeryrelationships}
\end{figure}

\begin{lemma}\label{lem:surgeryrelationships}
As shown in Figure~\ref{fig:surgeryrelationships}, $W(\alpha-1,\beta-1)=N(\alpha,\beta,+1)$.
\end{lemma}
\begin{proof}
This is a straightforward result of Kirby calculus.  See Chapters~4 and 5 of \cite{gompfstipsicz:4makc}, for example.
\end{proof}

\begin{lemma}[Martelli-Petronio,  \cite{mp:dfotm3m}]\label{lem:whiteheadlensspace}
$W(\alpha, \beta)$ is a lens space if and only if $\{\alpha,\beta\}$ is  $\{-1,-6+1/k\}$, $\{-2, -4+1/k\}$, $\{-3,-3+1/k\}$, or $\{p/q, \infty\}$ for some $k,p,q\in\Z$ and $(p,q)=1$. 
\end{lemma}
\begin{proof}
Using Lemma~\ref{lem:surgeryrelationships} this is a direct consequence of the results of Martelli and Petronio listed in Table~11 of  \cite{mp:dfotm3m}.
\end{proof}

\begin{lemma}[Proposition~3, \cite{hmw:sotwlygsm}]\label{lem:whiteheadfiber}
$W(p/q,\cdot)$ fibers over the circle if and only if $|q|\leq1$.  For each $p \in \Z$, $W(\pm p/1, \cdot)$ fibers with a once-punctured torus fiber and monodromy $\phi \cong \tau_x^{\mp p} \tau_y$.
\end{lemma}

\begin{remark}
If $q=0$ then $W(p/q, \cdot)$ is a solid torus and fibers over the circle with a disk fiber.  For $q \neq 0$, we refer the reader to \cite{hmw:sotwlygsm} for a proof.  Let us mention, nevertheless, that the monodromy is most directly viewed by regarding $W$ as the result of $-1$ surgery on one component of the Borromean rings.  Then $W(\pm p, \cdot)$ may be seen as the pair of surgeries $\pm p$ and $-1$ on the two-component chain link in the solid torus.
\end{remark}

\subsection{Heegaard splittings and height functions.}

A lens space $X$ has a unique genus one Heegaard splitting:  there is a torus, unique up to isotopy, that divides $X$ into two solid tori, e.g.\ \cite{bo:sdhdel}.  These are referred to as the {\em Heegaard torus} and {\em Heegaard solid tori} accordingly.
Endow $X$ with a height function $X \to [-\infty, +\infty]$ where the preimages of heights in $(-\infty, +\infty)$ are Heegaard tori and the preimages of $\pm \infty$ are the two unknotted circles.  Let $\hat{T}$ be a the preimage of $0$.  We refer to the Heegaard solid tori bounded by $\hat{T}$ and the objects in them as being {\em above} or {\em below} $\hat{T}$ depending on their relationship with respect to the height function on $X$.  The colors {\em amber} and {\em blue} will also be used respectively.  Denote the amber and blue Heegaard solid tori as $V_a$ and $V_b$.  Let $K_a = K \cap V_a$ and $K_b = K \cap V_b$.  Let $m_a$ and $m_b$ be meridians of $V_a$ and $V_b$ respectively.

\subsection{Grid number $1$ knots.}\label{sec:gridnumber}
Assume the meridians $m_a$ and $m_b$ of $V_a$ and $V_b$ respectively intersect minimally on $\hat{T}$.  Therefore in $L(r,q)$ they intersect in $r$ points.  The (unoriented) {\em grid number $1$ knots} are defined by unordered pairs of distinct points $\{z, w\}$ of $m_a \cap m_b$.  Let $D_a$ and $D_b$ be meridional disks of $V_a$ and $V_b$ with boundaries $m_a$ and $m_b$ respectively.  Given a pair of distinct points $\{z,w\}$ of $m_a \cap m_b$ there is a properly embedded arc $K_a$ and $K_b$ in $D_a$ and $D_b$ respectively with these points as their endpoints.  The knot $K=K_a \cup K_b$ is then a grid number $1$ knot.  

If in some direction along $m_b$ the two points $\{z,w\}$ are $n$ apart (i.e.\ there are $n-1$ points of $m_a \cap m_b$ in the interior of an arc of $m_b$ between $z$ and $w$) then we say $K$ is the {\em $n$th grid number $1$} knot in $L(r,q)$ along $m_b$.  Taken along $m_a$, $K$ is the $nq$th grid number $1$ knot where $nq$ is reduced modulo $r$.  Because we are working with knots in lens spaces up to homeomorphism we will understand that saying $K$ is the $n$th grid number $1$ knot means the count of $n$ is taken with respect to either $m_a$ or $m_b$.  In this article, for the torus knots at hand, the geometry will make apparent the correct choice; for the non-torus knots at hand, the symmetry makes the choice irrelevant.

\section{Torus knots}\label{sec:torusknots}

A torus knot $K$ in a lens space is an essential simple closed curve on $\hat{T}$.  The exterior of a torus knot is a Seifert fiber space over the disk with (at most) two exceptional fibers and is fibered over $S^1$ unless the curve is isotopic to a meridian of one of the Heegaard solid tori.  In this latter situation the curve bounds a disk in the lens space and is hence the trivial knot;  its exterior is the connect sum of the lens space and a solid torus.  

If $K$ is a torus knot in the lens space $X = V_a \cup_{\hat{T}} V_b$, then let $\delta_a = \Delta(K,m_a)$ and $\delta_b = \Delta(K,m_b)$ viewing the curves on the Heegaard torus.  Then we say that $K$ is a {\em type $\{\delta_a, \delta_b\}$ torus knot}.  The type of a torus knot reveals its element of homology with respect to the two canonical generators, the unknotted cores of $V_a$  and $V_b$.

\begin{lemma}\label{lem:torusknots}
If $K$ is a non-null homologous torus knot and an \opt--knot then either $K$ is unknotted or $E(K)$ fibers over the circle with a once-punctured torus fiber and, up to homeomorphism, for some integer $k\neq0$, either 
\begin{itemize}
\item[(I)] $K$ is a type $\{2,3\}$ torus knot in $L(6k-1,2k-1)$ arising as the core of the $-6+1/k$--surgery in $W(-1,-6+1/k)$ with monodromy $\phi \cong  \tau_x \tau_y$,
\item[(II)] $K$ is a type $\{2,4\}$ torus knot in $L(8k-2,4k+1)$ arising as the core of the $-4+1/k$--surgery in $W(-2,-4+1/k)$ with monodromy $\phi \cong  \tau_x^2 \tau_y$, or
\item[(III)] $K$ is a type $\{3,3\}$ torus knot in $L(9k-3,3k-2)$ arising as the core of the $-3+1/k$--surgery in $W(-3,-3+1/k)$ with monodromy $\phi \cong \tau_x^3 \tau_y$.
\end{itemize}
\end{lemma}

\begin{proof}
Since the exterior of a non-trivial (and hence non-null homologous) torus knot is fibered, if it is an \opt--knot with a compressible properly embedded once-puncutred torus in its exterior, then Lemma~\ref{lem:compress} says that the knot is unknotted.  Therefore let us now assume the once-punctured torus in the exterior of $K$ is incompressible.

The exterior $E(K)$ of a torus knot $K$ is the union of two solid tori $V'_a$ and $V'_b$ along essential annuli in their boundary.  (View $V_a$ and $V_b$ as the Heegaard solid tori above and below the Heegaard torus on which $K$ lies.  Dig out a trench from each along $K$ to form $V'_a$ and $V'_b$.)   Let $A =V'_a \cap V'_b$ be this annulus, and let $D_i$ be meridional disks of $V'_i$ minimally intersecting $A$, $i=a,b$.  It may be arranged so that $D_a$ and $D_b$ are disjoint.  Set $\delta_i = |\bdry D_i \cap A|$ and let $\lambda = \lcm(\delta_a, \delta_b)$.  We may now easily view $E(K)$ as a fiber bundle over the circle with fiber $S$ composed of $\lambda/\delta_1$ parallel copies of $D_a$ in $V'_a$ joined to $\lambda/\delta_b$ parallel copies of $D_b$ in $V'_b$ along $\lambda$ arcs in $A$.  In order for $S$ to be a once-punctured torus, we must have $\chi(S)=\lambda/\delta_a + \lambda/\delta_b - \lambda = -1$.  Rewritten as $1/\delta_a + 1/\delta_b +1/\lambda = 1$, the solutions are clear.  The only positive integer solutions are $\{\delta_a, \delta_b\} = \{2, 3\}$ and $\{\delta_a, \delta_b\} = \{2, 4\}$ if $\delta_a \neq \delta_b$ and $\delta_a=\delta_b=3$ otherwise.  Therefore the only possible torus knot types of \opt--knots are $\{2,3\}$, $\{2,4\}$, and $\{3,3\}$.  We now describe all such torus knots and their monodromies.  To facilitate this, regard the disks $D_i$ as regular ideal $\delta_i$--gons.

\begin{figure}
\centering
\input{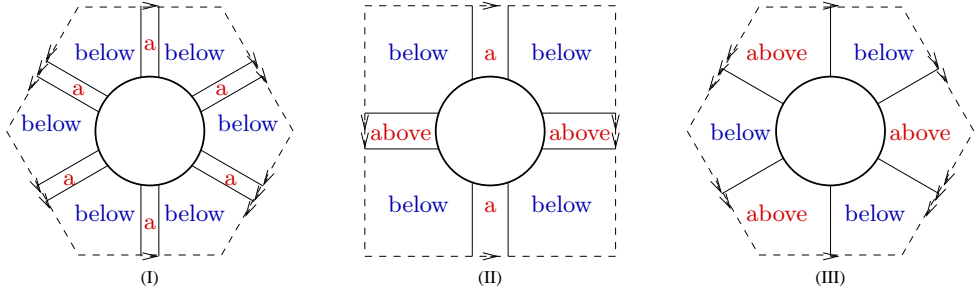}
\caption{The fibers of torus knots.}
\label{fig:torusknotfibers}
\end{figure}

Choosing $\delta_a=2$ and $\delta_b=3$, the fiber $S$ is assembled from $3$ copies of $D_a$ joined to $2$ copies of $D_b$ along $6$ arcs.  Sliding $S$ around $V'_a \cup V'_b$ through its monodromy back to itself permutes the $3$ copies of $D_a$ around $V'_a$ and the $2$ copies of $D_b$ around $V'_b$.  Doing so three times returns each copy of $D_a$ to itself with a $\pi$ rotation.  Doing so twice returns each copy of $D_b$ to itself with a $2\pi/3$ rotation.  Hence the monodromy may be now visualized as a $\pi/3$ rotation of a regular hexagon minus an open disk at the center with opposite sides identified; $\phi \cong \tau_x \tau_y$.  See Figure~\ref{fig:torusknotfibers}(I).  Hence by Lemma~\ref{lem:whiteheadfiber} $E(K) \cong W(-1, \cdot)$.  Martelli and Petronio then show in Table~11 of  \cite{mp:dfotm3m} that $W(-1,-6+1/k) \cong L(6k-1,2k-1)$ for $k\in\Z$ are the only lens space fillings of $W(-1,\cdot)$.

Choosing $\delta_a=2$ and $\delta_b=4$, the fiber $S$ is assembled from $2$ copies of $D_a$ joined to $1$ copy of $D_b$ along $4$ arcs.  The monodromy permutes the $2$ copies of $D_a$.  Done twice, each copy of $D_a$ returns to itself with a $\pi$ rotation.  Done once, $D_b$ returns to itself with a $\pi/2$ rotation.   Hence the monodromy may be visualized as a $\pi/2$ rotation of a square minus an open disk at the center with opposite sides identified; $\phi \cong \tau_x \tau_y \tau_x \cong \tau_x^2 \tau_y$.  See Figure~\ref{fig:torusknotfibers}(II).  Hence by Lemma~\ref{lem:whiteheadfiber} $E(K) \cong W(-2, \cdot)$.   Martelli and Petronio then show in Table~11 of \cite{mp:dfotm3m} that $W(-2,-4+1/k) \cong L(8k-2,4k+1)$ for $k\in\Z$ are the only lens space fillings of $W(-2,\cdot)$.

Choosing $\delta_a=3$ and $\delta_b=3$, the fiber $S$ is assembled from $1$ copy of $D_a$ joined to $1$ copy of $D_b$ joined along $3$ arcs.  The monodromy rotates each of the two disks by $2\pi/3$ (so that $|\bdry S| = 1$).  Hence the monodromy may be visualized as a $2\pi/3$ rotation of a regular hexagon minus an open disk at the center with opposite sides identified; $\phi \cong \tau_y \tau_x \tau_y \tau_x \cong \tau_x^3 \tau_y$.  See Figure~\ref{fig:torusknotfibers}(III).  Hence by Lemma~\ref{lem:whiteheadfiber} $E(K) \cong W(-3, \cdot)$.  Martelli and Petronio then show in Table~11 of  \cite{mp:dfotm3m} that for each $W(-3,-3+1/k) \cong L(9k-3,3k-2)$ for $k\in\Z$ are the only lens space fillings of $W(-3,\cdot)$.
\end{proof}

\begin{remark}
One may notice that the lens spaces $L(9k,3k+1)$ for $k\in \Z$ also contain type $\{3,3\}$ torus knots.  The fibers for these knots however are thrice-punctured spheres.
\end{remark}

\begin{lemma}\label{lem:torusknotgridnumber}
If $K$ is a knotted non-null homologous type $\{\delta_a, \delta_b\}$ torus knot and an \opt--knot then $K$ has grid number $1$ and is the $\delta_a$th or $\delta_b$th grid number $1$ knot.
\end{lemma}

\begin{proof}
Assume $K$ is a type $\{\delta_a, \delta_b\}$ torus knot.  
Let $m_a$ and $m_b$ be meridians of $V_a$ and $V_b$ that intersect minimally.  Then $m_a \cup m_b$ forms a toroidal grid.  Beginning at a vertex of the grid, travel along $m_b$ in one direction until the $\delta_a$th vertex (intersection of $m_a$) is reached.  
Since $K$ is a type $\{\delta_a, \delta_b\}$ torus knot, for one of the choices of directions along $m_b$, the original vertex must now be $\delta_b$ vertices (intersections of $m_b$) away in one direction along $m_a$.  Traveling in this direction back to the original vertex now completes a loop on $\hat{T}$ that is at least homotopic to $K$.  If the interior of the part of this loop that runs along $m_a$ is slightly lifted up into $V_a$ and the interior of the part that runs along $m_b$ is slightly pushed down into $V_b$, then the knot created is the $\delta_a$th grid number $1$ knot, counted along $m_b$.  Counting along $m_a$, this would be the $\delta_b$th grid number $1$ knot.

To show that this loop is actually simple and hence the torus knot we desire, we may simply show that no vertex of the grid in the loop is visited more than once.  Number the intersections of $m_a$ and $m_b$ from $0$ to $r-1$ in order along, say, $m_b$.  Then, since $X = L(r,q)$, the intersection points appear in order along $m_a$ as $0, q, 2q, \dots, (r-1)q$ reduced modulo $r$.  Therefore the loop on the grid encounters the following vertices in the order 
\[ 0,1, \dots, \delta_a,\delta_a + \dot{q}, \delta_a + 2 \dot{q}, \dots, \delta_a + \delta_b \dot{q}\]
 reduced modulo $r$ for some choice of $\dot{q} = \pm q^{\pm1}$.  
 
 Lemma~\ref{lem:torusknots} gives three families for which we must check that such a sequence exists with no repetitions (except at the beginning and end).  We do so to complete the proof.
\begin{itemize}
\item[(I)\quad]  
\begin{itemize}
\item[$k > 0$:] The type $\{2,3\}$ torus knot in $L(6k-1,2k-1)$ gives the sequence \[0,1,2, 2k+1, 4k, 6k-1.\]
\item[$k < 0$:] Set $k' = -k$.  The type $\{2,3\}$ torus knot in $L(6k'+1,2k'+1) \cong L(6k'+1,4k')$ gives the sequence \[0,1,2,4k'+2, 2k'+1,6k'+1.\]
\end{itemize}
\item[(II)\quad] 
\begin{itemize}
\item[$k > 0$:] The type $\{2,4\}$ torus knot in $L(8k-2,4k+1) \cong L(8k-2,2k-1)$ gives the sequence \[0,1,2,2k+1,4k, 6k-1,8k-2.\]
\item[$k < 0$:] Set $k' = -k$. The type $\{2,4\}$ torus knot in $L(8k'+2, 4k'-1)\cong L(8k'+2,6k'+1)$ gives the sequence \[0,1,2,6k'+3,4k'+2, 2k'+1, 8k'+2.\]
\end{itemize}
\item[(III)\quad]
\begin{itemize}
\item[$k>0$:] The type $\{3,3\}$ torus knot in $L(9k-3,3k-2)$ gives the sequence \[0,1,2,3,3k+1,6k-1,9k-3.\]
\item[$k<0$:] Set $k'=-k$.  The type $\{3,3\}$ torus knot in $L(9k'+3, 3k'+2)\cong L(9k'+3,6k'+1)$ gives the sequence \[0,1,2,3,6k'+4,3k'+2,9k'+3.\]
\end{itemize}
\end{itemize}
\end{proof}

\begin{remark}
If $K$ is unknotted in $L(r,q)$, then it may be identified as either the $1$st or the $q$th grid number $1$ knot in $L(r,q)$.
\end{remark}

\section{Fat-vertexed graphs of intersection}\label{sec:fatvertex}
Here we set up machinery used in the remainder of this article.   See Lickorish \cite{lickorish} for the introduction of fat-vertexed graphs in the study of intersecting surfaces.  See Gabai \cite{gabai:fatto3mIII}, Rieck \cite{rieck:hsomitdfs}, and Gordon's survey \cite{gordon:cmids} for the notion of thin position for knots in closed oriented $3$--manifolds and its interplay with intersecting surfaces.

Put $K$ in thin position with respect to the height function on our lens space $X$.  Arrange that $\hat{T}$ is a thick torus with respect to $K$ which $K$ intersects, say, $t$ times.  Let $H(K)$ be the closure of the tubular neighborhood $N(K)$.  Isotop $S$ rel--$\bdry$ in $E(K)$ to minimize $|S \cap T|$.  Abstractly cap off $S$ with a disk to form the torus $\hat{S}$.  Form the fat-vertexed graphs $G_S$ on $\hat{S}$ and $G_T$ on $\hat{T}$ as follows.  We will view our diagrams of $G_T$ from above.

The fat-vertex of $G_S$ is the disk abstractly capping off $S$.  The fat-vertices of $G_T$ are the disks $\hat{T} \cap H(K)$, the closures of the components of $\hat{T} \cap N(K)$.  We find it useful to regard the fat-vertex of $G_S$ as just $\bdry S$ and the fat-vertices as meridional disks of $H(K)$.  

Denote $T = \hat{T} - N(K) = \hat{T} \cap E(K)$.  The {\em edges} of $G_S$ and $G_T$ are the arcs of $S \cap T$.  Typically we do not include the closed curve components of $S \cap T$ in our graphs, but we will occasionally bring them into play.  The {\em faces} of $G_S$ and $G_T$ are the path closures of the components of the complements of the edges in $S$ and $T$ respectively.  Here too one may occasionally care to include the simple closed curves of $S \cap T$ with the edges in the formation of the faces.  A {\em corner} of a face $f$ is a component of $\bdry f \cap \bdry H(K)$.   An $n$--gon is a disk face with $n$ edges and $n$ corners.  For $n=1$, $2$, $3$,  and $4$ we call an $n$--gon a {\em monogon}, {\em bigon}, {\em trigon}, and {\em tetragon} respectively.  By the minimality assumption on $|S \cap T|$, if $K$ is not a torus knot,  then we may assume no face of $G_S$ or $G_T$ is a monogon (see e.g.\ Gordon~\cite{gordon:cmids}).  

Number the intersections of $K$ with $\hat{T}$ from $1$ to $t$ (taken modulo $t$) in order around $K$; the fat-vertices of $G_T$ are then labeled correspondingly.     Around $\bdry S$ each intersection with $\bdry T$ inherits one of these labels.     Observe that $\bdry S$ intersects each component of $\bdry T$ a total of $s$ times.  Consequentially, the labels $\{1, 2, \dots, t\}$ appear in order $s$ times around $\bdry S$.  Therefore $G_S$ and $G_T$ have $st/2$ edges.

Each edge $e$ of $G_T$ connecting vertices $i$ and $j$ is also an edge $e$ of $G_S$ with end points labeled $i$ and $j$.  We say $\{i,j\}$ is the {\em label pair} of the edge $e$ (as viewed in $G_S$).  Since $G_S$ has just one vertex, $\hat{T}$ is separating, and both $S$ and $T$ are orientable, the two labels of a label pair for an edge in $G_S$ must have opposite parity; all edges on $G_T$ connect vertices of opposite parity.   This is the Parity Rule, see \cite{cgls:dsok}, specialized to our context.  (Since $G_S$ has just one vertex we need not label it.   Edges in $G_T$ would all have the same label pair.  Edges in $G_S$ connect its single vertex to itself.)

Color each face $f$ of $G_S$ {\em amber} or {\em blue} according to whether a small collar neighborhood of $\bdry f$ in $f$ lies above or below $\hat{T}$ respectively.  This gives a ``checkerboard'' coloring of $G_S$.  We may also choose to include the simple closed curves of $S \cap T$ on $G_S$ and then color the resulting faces according to the color of the Heegaard solid torus in which their interiors lie.

Let $H_{i,i+1}$ be the $1$--handle of $H(K)$ running from vertex $i$ to vertex $i+1$ and $K_{i,i+1} = K \cap H_{i,i+1}$.  We say each $H_{i,i+1}$ and $K_{i,i+1}$ is {\em unknotted} if $K_{i,i+1}$ is contained in a meridional disk of the Heegaard solid torus in which it lies.  Denote $\bdry H_{i,i+1} \cap \bdry H(K)$ by $\bdry_v H_{i,i+1}$.  We will always assume that $H_{1,2}$ is amber.  Let $H_a = H(K) \cap V_a$ and $H_b = H(K) \cap V_b$.

A set of edges $\sigma=\{e_1, e_2, \dots, e_p\}$ in $G_S$ each with the same label pair that together bound a single disk face of $G_S$ is called a {\em Scharlemann cycle} of length $p$.   
The {\em label pair} of $\sigma$ is the label pair of the edges of $\sigma$.  

For any set of edges $\sigma$, if the subgraph of $G_T$ consisting of the edges of $\sigma$ and the vertices to which these edges are incident is contained in a disk on $\hat{T}$, then we say the edges of $\sigma$ {\em lie in a disk}.  If the subgraph is contained in an annulus on $\hat{T}$ but not in a disk, then we say the edges of $\sigma$ {\em lie in an essential annulus}.

We say two edges of $G_S$, $G_T$ are {\em parallel} if they cobound a bigon in $S$, $T$ respectively.  We say two faces $f_0$ and $f_1$ of $G_S$ are {\em parallel} if there is an isotopy of $f_0$ to $f_1$ in $E(K)$ keeping $\bdry f_0 \cap T$ on $T$ and $\bdry f_0 \cap \bdry H(K)$ on $\bdry H(K)$; in particular if $\Int (f_0 \cup f_1) \cap T = \emptyset$ then $f_0$ and $f_1$ cobound a product region in $E(K)$.

\section{Scharlemann cycles of length $2$ and $3$}\label{sec:schcycle}

\begin{lemma}[Lemma 2.1, \cite{gt:dsokwylsagok}]\label{lem:gt2.1}
Let $\sigma$ be a Scharlemann cycle in $G_S$ of length $p$ with label pair $\{i, i+1\}$ where $p$ is $2$ or $3$.  Let $f$ be the face of $G_S$ bounded by $\sigma$.  If the edges of $\sigma$ do not lie in a disk in $\hat{T}$ then they lie in an essential annulus $A$ in $\hat{T}$.  Furthermore, if $\Int f \cap \hat{T} = \emptyset$, then $M = N(A \cup H_{i,i+1} \cup f)$ is a solid torus such that the core of $A$ runs $p$ times in the longitudinal direction of $M$.
\end{lemma}

\begin{lemma}\label{lem:meridianandessannulus}
Continuing Lemma~\ref{lem:gt2.1} in the case that $\Int f \cap \hat{T} = \emptyset$, let $V$ be the solid torus bounded by $\hat{T}$ that contains $f$.  Then the core of $A$ runs $p$ times in the longitudinal direction of $V$.  In particular, the meridian of $V$ and the core of $A$ have distance $p$ on $\hat{T}$.
\end{lemma}

\begin{proof}
Since $A \cup H_{i,i+1} \cup f \subset V$ and $(A \cup H_{i,i+1} \cup f) \cap \hat{T} = A$, we may regard $M$ as its restriction to $V$ with $M \cap \hat{T} = A$.  Because the annulus $\bdry M - \Int A$ in $V$ is not parallel to $A$ through $M$, it must be parallel to $\hat{T} - \Int A$ through the solid torus $\overline{V - M}$.  Therefore a meridional disk of $M$ extends to a meridional disk of $V$.  Thus the core of $A$ will run $p$ times in the longitudinal direction of $V$ as well.  Hence the meridian of $V$ and the core of the annulus $A$ have distance $p$ on $\hat{T}$.
\end{proof}

\begin{lemma}\label{lem:gtS2cycles}
 Let $\sigma_i$ and $\sigma_j$ be two Scharlemann cycles of length 2 in $G_S$ with disjoint label pairs $\{i, i+1\}$ and $\{j, j+1\}$, respectively.  If the edges of each $\sigma_i$ and $\sigma_j$ lie in an essential annulus on $\hat{T}$, then $\sigma_i$ and $\sigma_j$ lie on opposite sides of $\hat{T}$.
\end{lemma}

\begin{proof}
Let $f_k$ be the face of $G_S$ bounded by $\sigma_k $, $k =i,j$.
Assume $f_i$ and $f_j$ have the same color, say blue.  Extend $f_i$ and $f_j$ across its corners radially to the cores of $H_{i,i+1}$ and $H_{j,j+1}$ to form the blue Mobius bands $g_i$ and $g_j$.  Since the edges of $f_k$ lie in an essential annulus on $\hat{T}$, so does $\bdry g_k$.   By disk exchanges, eliminate all simple closed curves of $\Int g_k \cap \hat{T}$ that bound disks in $\hat{T}$.  

Assume, for either $k=i$ or $j$, $\Int g_k \cap \hat{T}$  contains a simple closed curve that is essential on $\hat{T}$.  Then an innermost such curve on $g_k$ bounds a meridional disk $D$ in the amber solid torus; amber since $g_i$ and $g_j$ each intersect every blue meridional disk.  Because  $\bdry g_i$ and $\bdry g_j$ are both parallel to $\bdry D$, they bound disjoint amber meridional disks $D_i$ and $D_j$.  These disks may be taken to have interiors disjoint from $g_i$ and $g_j$.  Then $g_i \cup D_i$ and $g_j \cup D_j$ are disjoint embedded projective planes in a lens space.  This cannot occur.

Since $K$ intersects each $g_i$ and $g_j$ in a transverse blue arc, $\Int g_k \cap \hat{T}$ cannot contain a simple closed curve that is essential on $g_k$.  Therefore $\Int g_k \cap \hat{T} = \emptyset$.  Then $\bdry g_i$ and $\bdry g_j$ cobound an annulus $A$ on $\hat{T}$ so that together $g_i \cup A \cup g_j$ forms a Klein bottle embedded in the blue solid torus.  This cannot occur.  
\end{proof}

\begin{figure}
\centering
\input{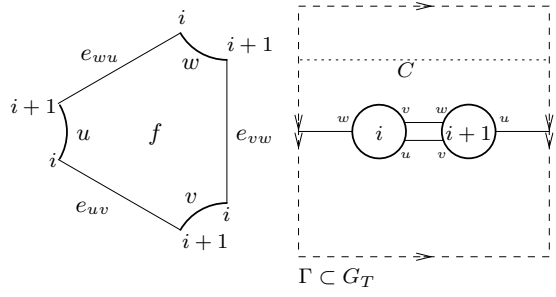}
\caption{The face of a Scharlemann cycle of length 3 with labeled corners and its edges on $\hat{T}$.}
\label{fig:meridianandS3}
\end{figure}

\begin{figure}
\centering
\input{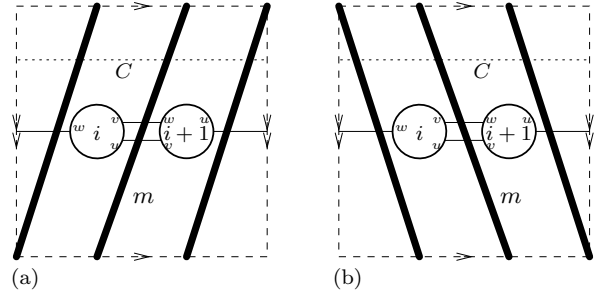}
\caption{A priori, two possible meridians for $V$ up to Dehn twists along $\gamma$.  The viable one is shown in Figure (a).}
\label{fig:twopossiblemeridians}
\end{figure}

\begin{lemma}\label{lem:meridiansandS3cycles}
Let $f$ be the face of a Scharlemann cycle of length 3 with label pair $\{i, i+1\}$ contained in the Heegaard solid torus $V$ whose edges do not lie in a disk on $\hat{T}$.  It induces the subgraph $\Gamma$ of $G_T$ on $\hat{T}$ as indicated in Figure~\ref{fig:meridianandS3} up to homeomorphism.   Assume $H_{i,i+1}$ is unknotted.   Then with respect to $\Gamma$ the meridian $m$ of a meridional disk of $V$ that is disjoint from $H_{i,i+1}$ appears on $\hat{T}$ as shown in Figure~\ref{fig:twopossiblemeridians}(a) (up to Dehn twists along an essential simple closed curve $C$ of $\hat{T}$ that is disjoint from $\Gamma$).
\end{lemma}


\begin{proof}
Figure~\ref{fig:meridianandS3} shows $f$ with a labeling of its corners and the induced labeling of the endpoints of the edges of $\Gamma$.
Let $c$ be an essential simple closed curve on $\hat{T}$ that is disjoint from $\Gamma$. By Lemma~\ref{lem:meridianandessannulus} $\Delta(m,C)=3$.  Therefore, up to Dehn twists along $c$, there are two possibilities for $m$.  These two possibilities are illustrated in Figure~\ref{fig:twopossiblemeridians}.

\begin{figure}
\centering
\input{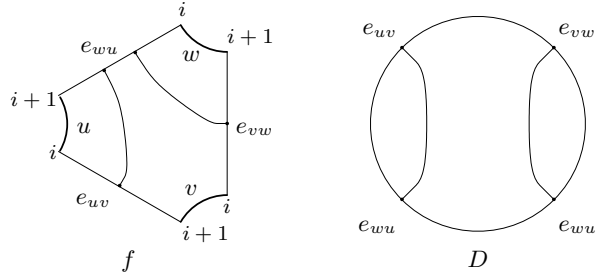}
\caption{The two arcs of intersection of the meridional disk $D$ with the face $f$.}
\label{fig:Dintersectsf}
\end{figure}

Let $D$ be a meridional disk of $V$ bounded by $m$.  By isotopies of $(D, \bdry D)$ in $(V, \bdry V)$ (thus keeping $m$ in $\hat{T}$), minimize $(|D \cap H_{i,i+1}|, |D \cap f|)$ lexicographically.   Since $H_{i,i+1}$ is unknotted, it is straightforward to check that once minimized  $| D \cap H_{i,i+1}|=0$ and $|D \cap f|=2$.  The two intersections of $D \cap f$ are arcs on $f$ connecting $e_{wu}$ to the other two edges, see Figure~\ref{fig:Dintersectsf}.  Hence  Figure~\ref{fig:twopossiblemeridians} also shows the two possibilities of how $m$ intersects $\Gamma$.   On $D$ these two arcs cut off disjoint subdisks.  Let $d$ be the one cut off by the arc on $f$ connecting $e_{wu}$ to $e_{vw}$.  Then $d \cap m$ is an arc on $\hat{T}$ that joins $e_{vw}$ to $e_{wu}$ at the point nearer vertex $i+1$ than the other intersection of $m$ with $e_{wu}$.  Therefore $m$ must lie as in Figure~\ref{fig:twopossiblemeridians}(a).
\end{proof}

\begin{lemma}\label{lem:2.2ish}
Let $\sigma$ be a Scharlemann cycle of length $p$ in $G_S$, and let $f$ be the face of $G_S$ bounded by $\sigma$.  If the edges of $\sigma$ lie in a disk then $r=p$.  Furthermore $\Int f \cap \hat{T} \neq \emptyset$, and there must be a component of $\Int f \cap \hat{T}$ that is essential on $\hat{T}$.
\end{lemma}

\begin{proof}
Let $\delta$ be a small disk in which the edges of $\sigma$ lie.  Let $\{i, i+1\}$ be the label pair of $\sigma$.  

Since lens spaces are irreducible, disk exchanges may be performed to replace $\hat{T}$ with another Heegaard torus $\hat{T}'$ such that $\Int f \cap \hat{T}'$ contains only curves which are essential on $\hat{T}'$.  Similarly $\delta$ may be replaced with a disk $\delta'$ in $\hat{T}'$ such that the edges of $\sigma$ lie in $\delta'$ and $\Int f \cap \delta' = \emptyset$.  Then $N(\delta' \cup H_{i,i+1} \cup f)$ forms a punctured lens space of order $p$.  Hence $r=p$.

If all components of $\Int f \cap \hat{T}$ were trivial on $\hat{T}$ then $\Int f \cap \hat{T}' = \emptyset$.  Then $N(\delta' \cup H_{i,i+1} \cup f)$ would form a punctured lens space within a neighborhood of one of the Heegaard solid tori.  This cannot occur.
\end{proof}

\begin{lemma}[cf.\ Lemma 2.3~\cite{gt:dsokwylsagok}]\label{lem:gt2.3}
Let $\sigma$ be a Scharlemann cycle of length $p$ in $G_S$ with label pair $\{i,i+1\}$ and let $f$ be the face of $G_S$ bounded by $\sigma$.  Suppose that $p \neq r$.  Then the edges of $\sigma$ cannot lie in a disk in $\hat{T}$.  Furthermore if $\Int f \cap \hat{T} \neq \emptyset$ then the components of $\Int f \cap \hat{T}$ are trivial curves in $\hat{T}$.
\end{lemma}

\begin{proof}
The proof of Lemma~2.3 of \cite{gt:dsokwylsagok} works for our purposes.  We have phrased our lemma to skirt Lemma~2.2 of \cite{gt:dsokwylsagok}.
Let us sketch the argument:

Assume $\sigma$ lies in a disk $D$ on $\hat{T}$.  By Lemma~\ref{lem:2.2ish} $p=r$, contrary to our hypothesis that $p \neq r$.  Therefore $\sigma$ cannot lie in a disk.

If $\Int f \cap \hat{T}$ contains a curve which is essential on $\hat{T}$ then $\sigma$ must lie in an essential annulus $A$.  Then, after suitable disk exchanges with $\hat{T}$, we may assume that there is a curve of $\Int f \cap \hat{T}$ that bounds a disk $\delta$ in $f$ which is a meridional disk of one of the solid tori bounded by $\hat{T}$.    We may then find an embedded annulus $A'$ which connects $\bdry \delta$ and a component of $\bdry A$ that is otherwise disjoint from $A \cup \delta \cup H_{i,i+1} \cup f$.  Then $N(A' \cup A \cup \delta \cup H_{i,i+1} \cup f)$ is a punctured lens space of order $p$.  This cannot occur within a lens space of order $r \neq p$.
\end{proof}

\begin{lemma}\label{lem:S2andS3}
Let $f_2$ and $f_3$ be disk faces of $G_S$ bounded by Scharlemann cycles $\sigma_2$ and $\sigma_3$ of length $2$ and $3$ respectively such that $f_2$ and $f_3$ have the same color and $\sigma_2$ and $\sigma_3$ have the same label pair $\{i,i+1\}$.  Let $V$ be the solid torus bounded by $\hat{T}$ of the same color as $f_2$ and $f_3$.  Let $m$ be its meridian.  Then
\begin{enumerate}
\item $\sigma_2$ and $\sigma_3$ each lie in an essential annulus $A_2$ and $A_3$ respectively on $\hat{T}$,
\item two edges of $\sigma_3$ and one edge of $\sigma_2$ are mutually parallel in the subgraph of $G_T$ consisting of the edges $\sigma_2 \cup \sigma_3$ and the vertices $i$ and $i+1$,
\item the cores $C_2$ and $C_3$ of the annuli $A_2$ and $A_3$ respectively have distance $1$ on $\hat{T}$, and 
\item the curves $m$, $C_2$, and $C_3$ may be oriented so that $[m] = 2[C_2] + 3[C_3] \in H_1(\hat{T})$.
\end{enumerate}
\end{lemma}

\begin{figure}
\centering 
\input{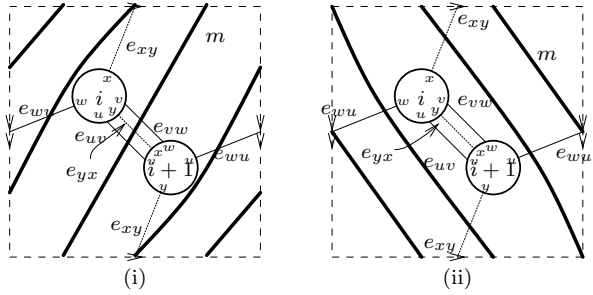}
\caption{Two possibilities for the meridian $m$ of a solid torus that contains the faces of Scharlemann cycles $\sigma_3 = \{e_{uv},e_{vw},e_{wu}\}$ and $\sigma_2 = \{e_{xy},e_{yx}\}$ on $G_T$.}
\label{fig:meridian}
\end{figure}

\begin{proof}
By Lemma~\ref{lem:gt2.3} for at least one of $p=2$ or $p=3$ the edges $\sigma_p$ cannot lie in a disk on $\hat{T}$.  Then by Lemma~\ref{lem:gt2.1} the edges of $\sigma_p$ must lie in an essential annulus $A_p$ on $\hat{T}$.  Let $p' \in \{2, 3\} - \{p\}$.  The proof of this lemma is now broken into the following three claims.

\begin{claim}\label{claim:1}
The edges of  $\sigma_{p'}$ lie in an essential annulus $A_{p'}$ on $\hat{T}$.
\end{claim}

\begin{proof} 
Assume the edges of $\sigma_{p'}$ lie in an disk.  
Then by Lemma~\ref{lem:2.2ish} both $r=p'$ and there must be a curve of $\Int f_{p'} \cap \hat{T}$ that is essential in $\hat{T}$.  Because $r \neq p$, Lemma~\ref{lem:gt2.3} implies that $\Int f_p \cap \hat{T}$ may only contain curves that are trivial in $\hat{T}$.  Perform disk exchanges to pass to a new Heegaard torus that is disjoint from $\Int f_{p}$ and intersects $\Int f_{p'}$ in only simple closed curves which are essential on the Heegaard torus.   We will retain the former names for all our objects that may have been affected by these disk exchanges.

By Lemma~\ref{lem:gt2.1} now $M=N(A_p \cup H_{i,i+1} \cup f_p)$ is a solid torus such that the core of $A_p$ runs $p$ times in the longitudinal direction of $M$.  
Also, now an innermost curve $\xi$ of $\Int f_{p'} \cap \hat{T}$ will bound a disk on $f_{p'}$ that is a meridional disk of one of the solid tori bounded by $\hat{T}$.
Since $\xi$ is essential on $\hat{T}$ it must be parallel to $A_p$.  Let $A'$ be an annulus on $\hat{T}$ giving a parallelism between $\xi$ and $A_p$.  Note that $A'$ is disjoint from $f_p$.  Then $N(\delta \cup A' \cup A_p \cup H_{x,x+1} \cup f_p)$ may be viewed as a $2$--handle $N(\delta \cup A')$ attached to $M$ along the core of $A_p$.  This produces a punctured lens space of order $p$ in our lens space of order $p'$.  This is a contradiction.
\end{proof}

We may identify $f_3$ and the subgraph of $G_T$ induced by $\sigma_3$ with $f$ and $\Gamma$ of Figure~\ref{fig:meridianandS3} so that $\sigma_3 = \{e_{uv},e_{vw},e_{wu}\}$.  Label the edges and corners of $f_2$ as shown in Figure~\ref{fig:labeledf2} so that $\sigma_2 = \{e_{xy},e_{yx}\}$.  
By Claim~\ref{claim:1} above, the edges $\sigma_3$ lie in the essential annulus $A_3$ on $\hat{T}$ and the edges $\sigma_2$ lie in the essential annulus $A_2$ on $\hat{T}$.  Let $C_2$ and $C_3$ be cores of $A_2$ and $A_3$ respectively.   

\begin{figure}
 \centering
 \input{labeledf2.pstex_t}
 \caption{}
 \label{fig:labeledf2}
\end{figure}

\begin{claim}\label{claim:2}  
Ignoring the curve $m$, the subgraph of $G_T$ induced by the edges $\sigma_2 \cup \sigma_3$ appears as in Figure~\ref{fig:meridian}.
Furthermore, $\Delta(C_2,C_3)=1$ on $\hat{T}$.
\end{claim}

\begin{proof}
Since $\sigma_2$ does not lie in a disk, we may assume $e_{xy}$, say, does not lie between $e_{uv}$ and $e_{vw}$ in $G_T$.  Therefore the endpoint of $e_{xy}$ on vertex $i$ either (1) appears between $e_{uv}$ and $e_{wu}$ (i.e.\ it is separated from the endpoint of $e_{vw}$ on vertex $i$ by the endpoints of $e_{uv}$ and $e_{wu}$)  or (2) between $e_{vw}$ and $e_{wu}$.  Since the corner $x$ of $f_2$ on $\bdry_v H_{i,i+1}$ is disjoint from the corners $u$,$v$, and $w$ of $f_3$, we may follow it along $\bdry_v H_{i,i+1}$ from the endpoint of $e_{xy}$ on vertex $i$ to the endpoint of $e_{yx}$ on vertex $i+1$ by following the corners of $f_3$ that flank it.

In case (1), the corner $x$ lies between the corners $w$ and $u$ on $\bdry_v H_{i,i+1}$ and so the endpoint of $e_{yx}$ on vertex $i+1$ must lie between $e_{vw}$ and $e_{wu}$.  Then in order for $e_{xy}$ and $e_{yx}$ to both have endpoints on the corner $y$ of $f_2$, this corner must also lie between the corners $w$ and $u$.  Consequentially, the cores $C_2$ and $C_3$ are isotopic on $\hat{T}$.  If $(\Int f_2 \cup \Int f_3) \cap \hat{T} \neq \emptyset$ then by Lemma~\ref{lem:gt2.3} the intersection may only consist of simple closed curves that are trivial on $\hat{T}$.  Eliminate these intersections by disk exchanges of $\hat{T}$ so that $(\Int f_2 \cup \Int f_3) \cap \hat{T} = \emptyset$.  Then by Lemma~\ref{lem:meridianandessannulus} for the meridian $m$ of the solid torus containing $f_2$ and $f_3$, $\Delta(m,C_2)=2$ and $\Delta(m,C_3)=3$.  This cannot occur since $C_2$ and $C_3$ are isotopic.

In case (2), the corner $x$ lies between the corners $w$ and $v$ on $\bdry_v H_{i,i+1}$ and so the endpoint of $e_{yx}$ on vertex $i+1$ must lie between $e_{uv}$ and $e_{vw}$.  Therefore the edge $e_{yx}$ lies in the bigon of $T$ bounded by $e_{uv}$ and $e_{vw}$.   The corner $y$ then lies between the corners $u$ and $v$ on $\bdry_v H_{i,i+1}$ so that the endpoint of $e_{xy}$ on vertex $i+1$ lies between $e_{uv}$ and $e_{wu}$.  Observe that $\Delta(C_2,C_3)=1$.  This completes the claim.
\end{proof}

\begin{claim}
The curves $m$, $C_2$, and $C_3$ may be oriented so that $[m] = 2[C_2] + 3[C_3] \in H_1(\hat{T})$.
\end{claim}

\begin{proof}
Because $\Delta(C_2,C_3)=1$ by Claim~\ref{claim:2}, $(\Int f_2 \cup \Int f_3) \cap \hat{T}$ may only consist of curves that are trivial on $\hat{T}$.  As in the proof of Claim~\ref{claim:1}, perform disk exchanges to pass to a new Heegaard torus that is disjoint from $\Int f_2 \cup \Int f_3$ and retain the former nomenclature.  Then Lemma~\ref{lem:meridianandessannulus} implies that $\Delta(C_2,m)=2$ and $\Delta(C_3,m)=3$.  Hence with respect to $C_2$ and  $C_3$ there are only two possibilities for the (unoriented) isotopy class of $m$ as described in the claim.
\end{proof}

These three claims complete the lemma.  The possible arrangements of the edges $\sigma_2=\{e_{xy},e_{yx}\}$ and $\sigma_3=\{e_{uv},e_{vw},e_{wu}\}$ with the meridian $m$ of $V$ are shown in Figure~\ref{fig:meridian}.
\end{proof}

\begin{figure}
\centering 
\input{isotopyontoT.pstex_t}
\caption{The isotopy rel--$\bdry$ of $K_{i,i+1}$ onto $\hat{T}$.}
\label{fig:isotopyontoT}
\end{figure}

\begin{lemma}\label{lem:ontoT}
Continuing Lemma~\ref{lem:S2andS3}, if we disregard the other arcs of $K \cap V$, the arc $K_{i,i+1}$ is isotopic rel--$\bdry$ in $V$ to an arc in $\hat{T}$ with interior disjoint from the edges $\sigma_2 \cup \sigma_3$ of $G_T$ that is incident to each vertex of $G_T$ between the consecutive edges of $\sigma_{3}$.  See Figure~\ref{fig:isotopyontoT}.
\end{lemma}

\begin{proof}
Among meridional disks of $V$ that are disjoint from $K_{x,x+1}$ assume $D$ is one that intersects $f_2 \cup f_3$ minimally.
Consider an arc $\gamma$ of $D \cap (f_2 \cup f_3)$ that is outermost on $D$ and the disk $\delta$ that it sections off.  By minimality, $\gamma$ cannot be boundary parallel on $f_2 \cup f_3$ to an edge of these faces.  Because the two edges $e_{uv}$ and $e_{vw}$ of $\sigma_3$ flank the edge $e_{yx}$ of $\sigma_2$ in $G_T$, $\gamma \in D \cap f_3$; $\gamma$ could not be outermost on $D$ otherwise.  Since $\gamma$ is not parallel to an edge of $f_3$ it must be parallel to a corner.  Join this rectangle of parallelism on $f_3$ to $\delta$ and perform a slight isotopy to create a bigon $B$ that has one edge on $\bdry_v H_{i,i+1}$, has its other edge on $\hat{T}$, and is otherwise disjoint from $H_{i,i+1} \cup f_2 \cup f_3$.  
Because $B$ is disjoint from $f_2 \cup f_3$ and two edges of $\sigma_3$ flank an edge of $\sigma_2$ in $G_T$, the edge of $B$ on $\bdry_v H_{i,i+1}$ must lie between the two consecutive corners of $f_3$.  Thus $K_{i,i+1}$ is isotopic rel--$\bdry$ to the edge of $B$ on $\hat{T}$ as shown in Figure~\ref{fig:isotopyontoT}.
\end{proof}

\begin{lemma}\label{lem:twoS3cyclesinsolidtorus}
Let $f$ and $f'$ be the faces of two order 3 Scharlemann cycles in $G_S$ of the same color and with the same label pair $\{i, i+1\}$.  Assume $\Int (f \cup f') \cap \hat{T} = \emptyset$.  Let $V$ be the solid torus that contains $f$ and $f'$.  Then $f$ and $f'$ are parallel in $V-\Int H_{i, i+1}$.
\end{lemma}

\begin{proof}
Since $M=V-(N(f) \cup \Int H_{i,i+1})$ is a solid tours on which $\hat{T} \cap \bdry M$ is the disjoint union of a disk $D$ and an annulus $A$ that is longitudinal.  Furthermore $\bdry_v H_{i,i+1} \cap \bdry M$ is three rectangles, one connecting the two components of $\bdry A$ and the other two connecting the two components of $\bdry A$ to $\bdry D$.  Since $f'$ is a properly embedded disk in $M$ (with boundary in $(\hat{T} \cup \bdry_v H_{i,i+1}) \cap \bdry M$), it must either be boundary parallel or a meridional disk of $M$.  

Since $\bdry f'$ traverses the rectangles of $\bdry_v H_{i,i+1} \cap \bdry M$ a total of three times it either traverses the three rectangles each once or it traverses the rectangle connecting the two components of $\bdry A$ thrice. But since one of these rectangles may be traversed by $\bdry f'$ only in a single direction, if the latter occurs then $\bdry f'$ intersects $\bdry A$ in three spanning arcs; no disk in a solid torus has such a boundary.  Therefore the former case occurs and $\bdry f'$ must intersect $A$ twice (and $D$ once).  Therefore $f'$ cannot be a meridional disk.  Hence $f'$ is boundary parallel, and it thus follows that $f'$ is parallel to $f$. 
\end{proof}

\begin{lemma}\label{lem:disjtinterior}
Assume $t=2$.  If $s > 2$ then the interior of a bigon of $G_S$ is disjoint from $\hat{T}$.  If $s >3$ then the interior of a trigon of $G_S$ is disjoint from $\hat{T}$.
\end{lemma}

\begin{proof}
Let $\sigma$ be two or three edges of $G_S$ that bound a $p$--gon face $f$ with $\Int f \cap \hat{T} \neq \emptyset$, $p=2$ or $3$ respectively.  Because $t=2$, $\sigma$ is a Scharlemann cycle of length $p$.  Since $r\geq s> p$, by Lemma~\ref{lem:gt2.3} the edges of $\sigma$ do not lie in a disk.   By Lemma~\ref{lem:gt2.1} the edges of $\sigma$ lie in an essential annulus $A$ on $\hat{T}$.  Furthermore by Lemma~\ref{lem:gt2.3} the components of $\Int f \cap \hat{T}$ are trivial curves in $\hat{T}$.   Indeed, since $t=2$ these components must be trivial curves in $T$; otherwise $G_T$ would lie in a disk.  Since lens spaces are irreducible and $K$ is not contained in a $3$--ball, we may isotop $\Int f$ to remove its intersections with $\hat{T}$.  This reduces $|S \cap T|$ contrary to our minimality assumptions.
\end{proof}

\section{$s\geq3$ and $t=2$}\label{sec:s=3t=2}
This section is devoted to proving the following theorem.
\begin{thm}\label{thm:sgeq3}
If $s \geq 3$ then $s=3$ and, for each integer $k \neq 0$, $K$ is a knot in $W(-3,-3+1/k) \cong L(9k-3,3k-2)$ arising as the core of the $-3$ surgery.  
Equivalently, up to homeomorphism, $K$ is the $|3k-1|$th grid number $1$ knot in $L(9k-3,3k-2)$.  
Furthermore $E(K)$ is fibered if and only if $k= \pm1$:
\begin{itemize}
\item   If $k=+1$ then $K$ is a torus knot in $L(6,1)$ with monodomy $\phi \cong \tau_x^2 \tau_y$.
\item  If $k=-1$ then $K$ is a non-torus knot in $L(12,5)$ monodromy $\phi \cong \tau_x^4 \tau_y$.
\end{itemize}
\end{thm}
\begin{proof}
We give an overview of the proof.  Assume $s \geq 3$.  Theorem~1.1 of \cite{baker:sgkilshsbn} states that a knot in a lens space with an essential once-punctured genus $g$ surface properly embedded in its exterior such that the boundary of the surface is distance at least $4g-1$ from the meridian (i.e.\ $4g-1 \leq s$) then the knot is $1$--bridge with respect to the genus one Heegaard splitting of the lens space.  Since our surface $S$ has genus $1$ and we are assuming $s \geq 3$, it follows that $K$ is $1$--bridge.  Therefore $t=2$.  

With $t=2$, the first three subsections here conclude that $s=3$ and that $G_S$ must appear as in Figure~\ref{fig:feweredges}(b).  Then Lemma~\ref{lem:somegenusoneknots} describes the lens spaces which contain knots with such surfaces $S$.  
Lemma~\ref{lem:order3gridnumber} identifies each knot $K$ as grid number $1$.  Lemma~\ref{lem:order3surgdesc} gives the surgery description of this family of knots.  Lemma~\ref{lem:somegenusonefibered} shows that $E(K)$ is fibered if and only if $k = \pm1$ and that $S$ is an incompressible if and only if $k \neq 0$.  Lemma~\ref{lem:order3fiber} shows that if $k=+1$ then $K$ is a torus knot.  Lemma~\ref{lem:order3kplus1} shows that if $k=-1$ then $K$ has fibered exterior with monodromy $\phi \cong \tau_x^4 \tau^y$.
\end{proof}

Throughout this section we assume $s \geq 3$.  As stated in the proof above for use in the ensuing lemmas, $K$ is therefore $1$--bridge.  Hence $K=K_a \cup K_b$ where $K_a$ and $K_b$ are unknotted arcs in $V_a$ and $V_b$ respectively.  Consequentially, $t=2$.

\begin{figure}
\centering
\input{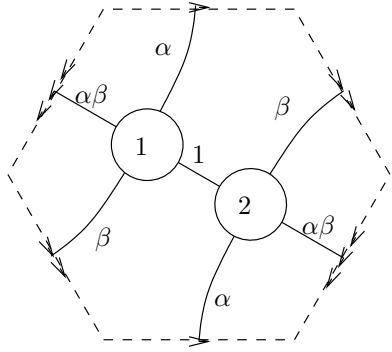}
\caption{The four edge class labels of $G_T$.}
\label{fig:4edgeclasslabels}
\end{figure}

Since $t=2$, the graph $G_T$ has two vertices and $s$ edges.  By the Parity Rule, there are at most four edge classes in $G_T$ which we denote $1, \alpha, \beta, \alpha\beta$ (see \S4 \cite{gt:dsokwylsagok} and \cite{gl:dsokcetI}) as shown in Figure~\ref{fig:4edgeclasslabels}.  Label an edge of $G_S$ by the class of the corresponding edge in $G_T$.  The label of an edge is called its {\em edge class label}.

\subsection{Bigons of $G_S$}

The following three lemmas give strong restrictions on the graph $G_S$.

\begin{lemma}[Lemma 4.1, \cite{gt:dsokwylsagok}]\label{lem:gt4.1}
Since $s > 2$, any two bigons in $G_S$ of the same color have the same pair of edge class labels.  Furthermore, the two bigons are parallel.
\end{lemma} 

\begin{proof}
Since the interior of a bigon is disjoint from $\hat{T}$ by Lemma~\ref{lem:disjtinterior}, the proof of Lemma~5.2 \cite{gl:dsokcetI} continues to hold for us too.  Let us sketch the argument.

Let $B_1$ and $B_2$ be two bigons of $G_S$ in, say, $V_a$.  Recall $H_{1,2}$ is the $1$--handle neighborhood of $K$ in $V_a$.  Then $M=V_a - (N(B_1) \cup \Int H_{1,2})$ is a solid torus on which $\hat{T} \cap \bdry M$ is a longitudinal annulus.  Since $B_2$ is a properly embedded disk in $M$, it must either be boundary parallel or a meridional disk.  A meridional disk will meet $\bdry M - \hat{T}$ minimally once.  Because the corners of $B_2$ must cross $\bdry_v H_{1,2}$ twice, $B_2$ must be boundary parallel.  It follows that $B_2$ is indeed parallel to $B_1$.  Hence it has the same pair of edge class labels.
\end{proof}

\begin{lemma}[cf.\ Lemma 4.2, \cite{gt:dsokwylsagok}] \label{lem:gt4.2}
Parallel edges of $G_S$ have distinct edge class labels.
\end{lemma}

\begin{proof}
Assume two parallel edges $e_1$ and $e_2$ of $G_S$ have the same edge class labels with no third edge in between also having the same edge class label.  By Lemma~\ref{lem:gt2.3} these two edges cannot be adjacent in $G_S$.  Because Lemma~\ref{lem:gt4.1} implies that two bigons of the same color must have the same pair of edge class labels, there may be at most two edges between $e_1$ and $e_2$ in $G_S$.  If there are these two edges between $e_1$ and $e_2$, then they must have the same edge class label.  This however contradicts Lemma~\ref{lem:gt2.3}.  Therefore there is necessarily just one edge $e_{1\frac{1}{2}}$ of $G_S$ between $e_1$ and $e_2$, and this edge has a different edge class label than that of $e_1$ and $e_2$.

Let $B_i$ be the bigon of $G_S$ bounded by $e_i$ and $e_{1\frac{1}{2}}$.  
 Extend $B_1$ and $B_2$ radially to the core $K$ of $H_{1,2} \cup H_{2,1}$ to form two M\"obius bands $A_1$ and $A_2$ respectively on either side of the torus $\hat{T}$.  Since $e_1$ and $e_2$ are in the same edge class, $\bdry A_1 \cup \bdry A_2$ bounds a bigon $B \subset \hat{T}$.
Thus $A_1 \cup A_2 \cup B$ is a Klein bottle $R$ in our lens space $X$ on which $K$ lies.  Moreover $K$ lies as an essential curve on $R$ and $R-N(K)$ is an annulus.  Since the curves of $\bdry (R-N(K))$ have induced coherent orientations on $\bdry N(K)$, then $2[K]=0$ in $H_1(X)$.  Hence $s=2$ contrary to our assumption that $s \geq 3$.
\end{proof}

\begin{remark}
Lemma 4.2 of \cite{gt:dsokwylsagok} appeals to a hypothesis that $K$ is hyperbolic for its proof.  Our proof uses that $s\geq3$ instead.
\end{remark}

\begin{lemma}[Lemma 4.3,  \cite{gt:dsokwylsagok}] \label{lem:gt4.3}
Since $s>2$, $G_S$ cannot contain more than $3$ mutually parallel edges.
\end{lemma}

\begin{proof}
This is an easy consequence of Lemmas~\ref{lem:gt4.1} and \ref{lem:gt4.2}.
\end{proof}

\subsection{If $G_S$ has $3$ mutually parallel edges...}
Since $S$ is a once-punctured torus, $G_S$ may have at most three edge classes.  By Lemma~\ref{lem:gt4.3} there are at most $3$ edges in any $G_S$--edge class.  In this subsection we assume at least one $G_S$--edge class does contain $3$ mutually parallel edges.

\begin{lemma}\label{lem:no3paralleledges}
If $s \geq 3$ and $K$ is not a torus knot then $G_S$ cannot have $3$ mutually parallel edges.
\end{lemma}

\begin{proof}
After enumerating all possibilities for $G_S$ having at most $3$ edges in any $G_S$--edge class (and $3$ edges in at least one $G_S$--edge class) and eliminating those that violate the Parity Rule, there are four left.  These possibilities are shown in Figure~\ref{fig:threeparallel}.  

\begin{figure}
\centering
\input{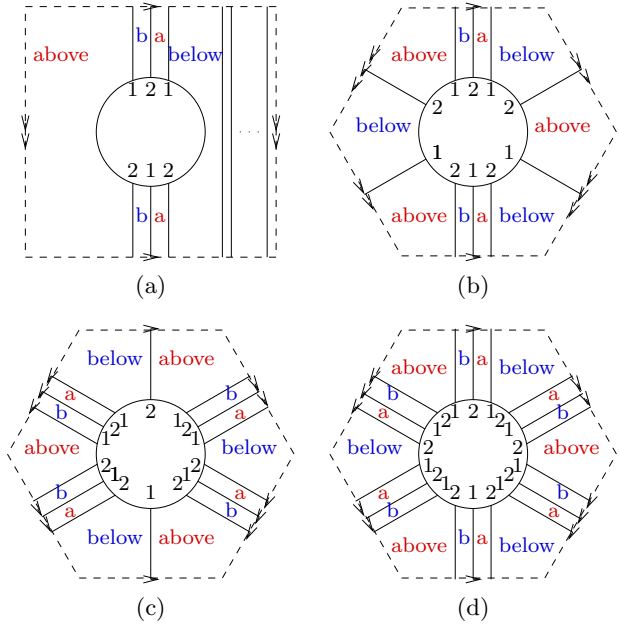}
\caption{If $G_S$ has three mutually parallel edges, no set of more than three mutually parallel edges, and obeys the Parity Rule then, a priori, $G_S$ may be one of the four shown.}
\label{fig:threeparallel}
\end{figure}

We eliminate the remaining possibilities with Claim~\ref{claim:3pa}, Claim~\ref{claim:3pc}, Claim~\ref{claim:3pd}, and Lemma~\ref{lem:3pb} below.
\end{proof}

\begin{claim}\label{claim:3pa}
Figure~\ref{fig:threeparallel}(a) cannot occur.
\end{claim}

\begin{proof}
In Figure~\ref{fig:threeparallel} (a), there must be an odd (and hence nonzero) number of essential simple closed curves of $S \cap T$ on $G_S$ in order to achieve the checkerboard coloring.  Such simple closed curves must also be essential in $T$ since $G_T$ does not lie in a disk and $S$ is incompressible.  However, since Lemma~\ref{lem:gt4.2} states that the three parallel edges of $G_S$ must each have different $G_T$--edge class labels, the complement of these three arcs on $T$ must be a disk.  Hence there can be no essential simple closed curve of $S \cap T$ on $T$.  This is a contradiction.  Figure~\ref{fig:threeparallel} (a) cannot occur.
\end{proof}

The remaining three possibilities for $G_S$ each contain for each color a trigon and one, two, or three bigons.  By Lemma~\ref{lem:gt4.2} the three $G_T$--edge class labels of three mutually parallel edges of $G_S$ must be distinct.  Then by Lemma~\ref{lem:gt4.1} the $G_T$--edge class labels for one set of three mutually parallel edges determines the $G_T$--edges class labels for every other set of three mutually parallel edges.  Furthermore, by Lemma~\ref{lem:S2andS3} a trigon of $G_S$ must have two edges in the same $G_T$--edge class as an edge of a bigon of $G_S$ of the same color.  This then implies that the remaining two edges of the trigon and bigon are labeled with two other distinct $G_T$--edge class labels.

\begin{claim}\label{claim:3pc}
Figure~\ref{fig:threeparallel}(c) cannot occur.
\end{claim}

\begin{proof}
After labeling the sets of three mutually parallel edges in accordance with Lemma~\ref{lem:gt4.2} and Lemma~\ref{lem:gt4.1}, there is no choice of label for the final edge that is consistent with Lemma~\ref{lem:S2andS3}.
\end{proof}

\begin{claim}\label{claim:3pd}
Figure~\ref{fig:threeparallel}(d) cannot occur.
\end{claim}

\begin{proof}
The pattern of $G_T$--edge class labels forced upon Figure~\ref{fig:threeparallel}(d) by the sets of three mutually parallel edges causes the two trigons to have all three edges with the same $G_T$--edge class label.  Because in this case $r \geq s=9 $, this contradicts Lemma~\ref{lem:gt2.3}.
\end{proof}

\begin{figure}
\centering
\input{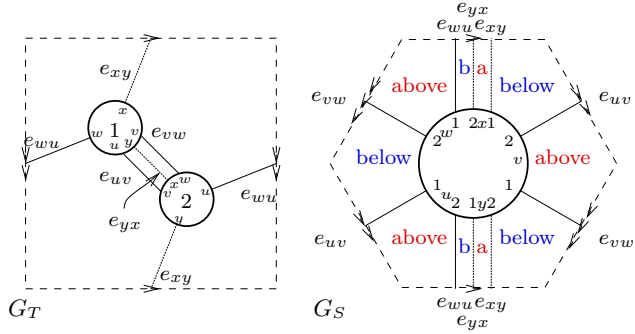}
\caption{Labels for Figure~\ref{fig:threeparallel}(b).}
\label{fig:threeparallelblabel}
\end{figure}

\begin{lemma}\label{lem:3pb}
If $G_S$ appears as in Figure~\ref{fig:threeparallel}(b), then $K$ is a torus knot.
\end{lemma}

\begin{proof}
Assume $G_S$ does appear as in Figure~\ref{fig:threeparallel}(b).
There are two Scharlemann cycles of length $3$, $\sigma_{3a}$ and $\sigma_{3b}$, bounding amber and blue faces $f_{3a}$ and $f_{3b}$ respectively; there are two Scharlemann cycles of length $2$, $\sigma_{2a}$ and $\sigma_{2b}$, bounding amber and blue faces $f_{2a}$ and $f_{2b}$ respectively.  

By Lemma~\ref{lem:S2andS3} we may assume $\sigma_{3a} = \{e_{uv}, e_{vw}, e_{wu}\}$ and $\sigma_{2a} = \{e_{xy} e_{yx}\}$ labels the edges as in Figure~\ref{fig:meridian}.  As there are only these five edges in $S \cap T$, Lemma~\ref{lem:S2andS3} also implies that $\sigma_{3b} = \{e_{uv}, e_{vw}, e_{xy}\}$ and $\sigma_{2b} = \{e_{wu}, e_{yx}\}$.  With this labeling, $G_T$ and $G_S$ appear as in Figure~\ref{fig:threeparallelblabel}.

Lemma~\ref{lem:S2andS3} describes how the meridians of $V_a$ and $V_b$ may appear on $\hat{T}$ with respect to $G_T$.
The amber meridian $m_a$ appears as $m$ in Figure~\ref{fig:meridian}(a) or (b).  Similarly, the blue meridian $m_b$ appears on $\hat{T}$ with respect to $G_T$ as a reflection of $m$ across a line of slope $1$ in Figure~\ref{fig:meridian}(a) or (b).   

Put together, the four possibilities for $m_a$ and $m_b$ on $\hat{T}$ are shown in Figure~\ref{fig:twomeridians}.  The latter two are reflections of the former two.  In (ii) and (iii), $m_a$ and $m_b$ define a lens space of order $r=13$ which cannot contain our knot with $s=5$.  The meridians $m_a$ and $m_b$ in (i) and (iv) each define the lens space $L(5,1)$ which may contain our knot.  These meridians (in relation to $G_T$)  are both of the type shown in Figure~\ref{fig:meridian}(i).

\begin{figure}
\centering
\input{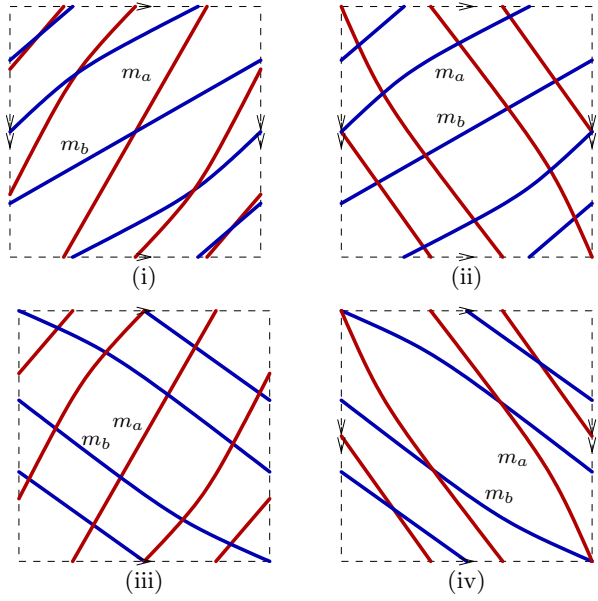}
\caption{The four possibilities for the two meridians of the Heegaard solid tori for Figure~\ref{fig:threeparallelblabel}.}
\label{fig:twomeridians}
\end{figure}

\begin{figure}
\centering
\input{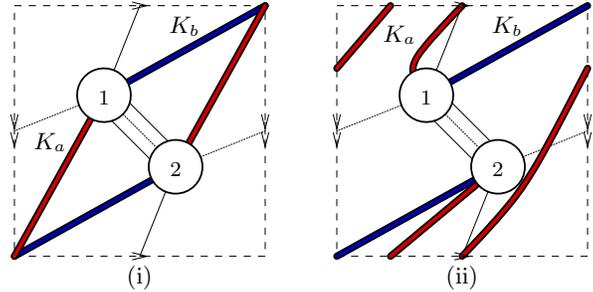}
\caption{(i) Isotopies of $K_a$ and $K_b$ onto $\hat{T}$.  (ii) Result of isotopy of $K_a$ across a meridional disk of $V_a$ back onto $\hat{T}$.  $K$ is thus a torus knot.}
\label{fig:ontoT}
\end{figure}

Lemma~\ref{lem:ontoT} allows us to obtain a diagram of $K$ on $\hat{T}$ in relationship to $G_T$.  The amber arc $K_a$ of $K$ is isotopic rel--$\bdry$ in $V_a$ to an arc in $\hat{T}$ with interior disjoint from the edges of $G_T$ that is incident to each vertex of $G_T$ between the consecutive edges of $\sigma_{3a}$.  Similarly, the blue arc $K_b$ of $K$ is isotopic rel--$\bdry$ in $V_b$ to an arc in $\hat{T}$ with interior disjoint from the edges of $G_T$ that is incident to each vertex of $G_T$ between the consecutive edges of $\sigma_{3b}$.    The result of these isotopies is shown in Figure~\ref{fig:ontoT}(i).  Note that $K_a$ lies above $K_b$.  The arc $K_b$, say, now may be isotoped rel--$\bdry$ across the meridional disk of $V_b$ so that it lies on $\hat{T}$ overlapping $K_a$ only at their endpoints as shown in Figure~\ref{fig:ontoT}(ii).  Hence $K$ is a torus knot.
\end{proof}

\begin{remark}
In fact if $G_S$ does appear as in Figure~\ref{fig:threeparallel}(b), then $K$ is a type $\{2,3\}$ torus knot in $L(5,1)$.
\end{remark}

\subsection{If $G_S$ only has fewer than $3$ mutually parallel edges...}
We now may assume no edge class of $G_S$ has more than $2$ mutually parallel edges.  

\begin{figure}
\input{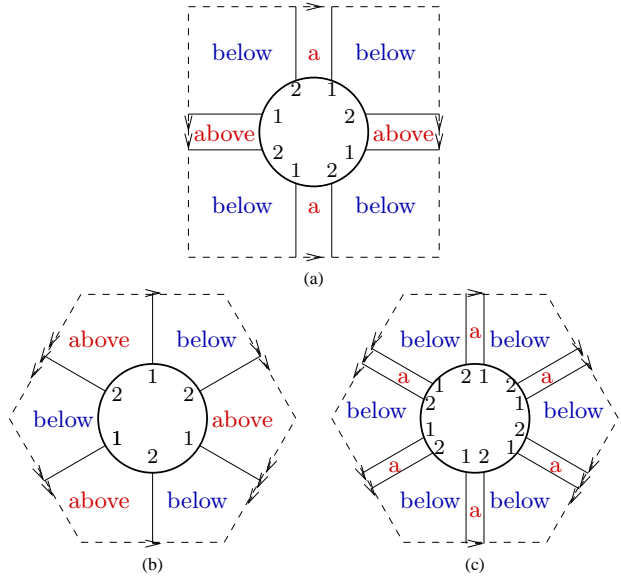}
\caption{These are the three possibilities for $G_S$ that have no set of three mutually parallel edges, obey the parity rule, and have $s \geq 3$.}
\label{fig:feweredges}
\end{figure}

\begin{lemma}\label{lem:no2paralleledges}
If $s \geq 3$ and $G_S$ does not have $3$ mutually parallel edges then $G_S$ appears as in Figure~\ref{fig:feweredges}(b).
\end{lemma}

\begin{figure}
\input{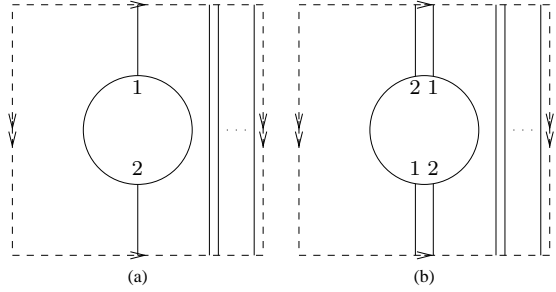}
\caption{These two possibilities for $G_S$ have no set of three mutually parallel edges and obey the parity rule, yet have $s < 3$.}
\label{fig:toofeweredges}
\end{figure}

\begin{proof}
After enumerating all possibilities for $G_S$ and eliminating those that violate the parity rule, there are five left.  The two shown in Figure~\ref{fig:toofeweredges} correspond to $s=1$ and $s=2$, respectively, contrary to our assumption that $s \geq 3$.  Therefore we are left to consider the three remaining possibilities shown in Figure~\ref{fig:feweredges}.  We eliminate Figure~\ref{fig:feweredges}(c) in Lemma~\ref{lem:2pc} and Figure~\ref{fig:feweredges}(a) in Lemma~\ref{lem:2pa} below leaving Figure~\ref{fig:feweredges}(b) as the only viable possibility.  
\end{proof}
We determine the knots $K$ with $G_S$ as in Figure~\ref{fig:feweredges}(b) in \S\ref{sec:2pb}.

\begin{lemma}\label{lem:2pc}
$G_S$ cannot appear as in Figure~\ref{fig:feweredges}(c).
\end{lemma}

\begin{proof}
Assume $G_S$ does appear as in Figure~\ref{fig:feweredges}(c).  The graph $G_S$ contains two Scharlemann cycles of length $3$, say $\sigma_3$ and $\sigma_3'$, of the same color and with the same label pair.  Since $s=6$, by Lemma~\ref{lem:disjtinterior} the faces of $\sigma_3$ and $\sigma_3'$ must have interior disjoint from $\hat{T}$.  Then by Lemma~\ref{lem:twoS3cyclesinsolidtorus} the faces of $\sigma_3$ and $\sigma_3'$ are parallel.  Hence the six edges of $G_T$ lie in an essential annulus on $\hat{T}$; four in one $G_T$--edge class and two in another.

The graph $G_S$ also contains three Scharlemann cycles of length $2$, say $\sigma_2$, $\sigma_2'$, and $\sigma_2''$, of the other color.  By Lemma~\ref{lem:gt4.1} all three have the same pair of $G_T$--edge class labels.  Hence three edges of $G_T$ lie in one $G_T$--edge class and the other three lie in another.  This is a contradiction.
\end{proof}

\begin{lemma}\label{lem:2pa}
$G_S$ cannot appear as in Figure~\ref{fig:feweredges}(a).
\end{lemma}

\begin{proof}
Assume $G_S$ does appear as in Figure~\ref{fig:feweredges}(a).  Then $G_S$ contains two Scharlemann cycles of length $2$, say $\sigma_2$ and $\sigma_2'$, with faces $f_2$ and $f_2'$ respectively, say colored amber, and the same label pair.   By Lemma~\ref{lem:gt4.1} the faces of $\sigma_2$ and $\sigma_2'$ must be parallel.  Let $\sigma_2 = \{e_{xy}, e_{yx}\}$ with  corners $x$ and $y$ as in Figure~\ref{fig:labeledf2} (setting $i=1$).  Similarly let $\sigma_2' = \{e_{uv}, e_{vu}\}$ with corners $u$ and $v$ -- replace $x$ with $u$ and $y$ with $v$ in Figure~\ref{fig:labeledf2}.  Without loss of generality we may assume these edges and corners appear on $G_S$ and $G_T$ as shown in Figure~\ref{fig:fe2(c)}.  The two pairs of parallel edges in $G_T$ hands us two bigons $D_1$ and $D_2$ on $T$.  

\begin{figure}
\centering
\input{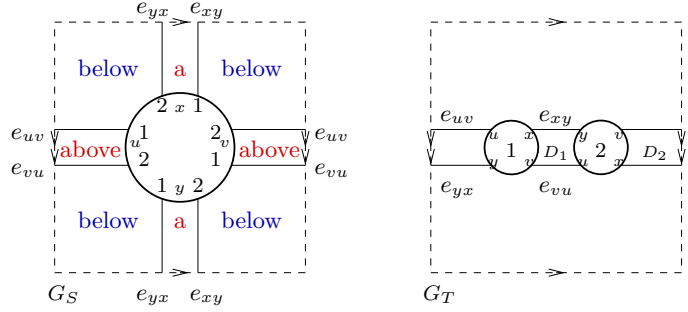}
\caption{The edges of Figure~\ref{fig:feweredges}(a) on $G_S$ and $G_T$ with labelings.}
\label{fig:fe2(c)}
\end{figure}

The graph $G_S$ also contains a Scharlemann cycle $\sigma_4$ of length $4$ bounding a blue face $f_4$.  Its four edges coincide with those of the two order two Scharlemann cycles.  

A simple closed curve $\gamma$ of $\Int f_4 \cap \hat{T} \neq \emptyset$ that is innermost on $f_4$ bounds a disk in $f_4$ and must also bound a disk in $T$.  Otherwise $\gamma$ would be parallel on $\hat{T}$ to the core of the annulus in which the edges of $\sigma_2$ lie.  By Lemma~\ref{lem:gt2.1} $\gamma$ runs twice in the longitudinal direction of $V_a$, hence the disk in $f_4$ that $\gamma$ bounds lies in $V_b$.  A neighborhood of this disk together with $V_a$ forms a punctured lens space of order $2$ embedded in our lens space $X$.  Hence $s=2$ contrary to our assumptions.
Since $\gamma$ bounds a disk in $T$ an isotopy of $\Int S$ will reduce $|S \cap T|$.  Hence $\Int f_4 \cap \hat{T} = \emptyset$.

We may now form an embedded once-punctured Klein bottle $R'$ in $V_b$ by attaching $D_1$ and $D_2$ to $f_4$ along their common edges.  Since $\bdry R' \subset \bdry_v H_{2,1}$ it must either be a meridian of $K$ or bound a disk on $\bdry_v H_{2,1}$.  In either case there is a properly embedded disk in $H_{2,1}$ with the same boundary as $R'$.  Together, $R'$ and this disk form an embedded Klein bottle in the solid torus $V_b$.  This cannot occur.
\end{proof}

\subsection{When $G_S$ appears as Figure~\ref{fig:feweredges}(b).}\label{sec:2pb}

\begin{lemma} \label{lem:notinadisk}
If $G_S$ is as in Figure~\ref{fig:feweredges}(b) then $G_T$ does not lie in a disk on $\hat{T}$.
\end{lemma}

\begin{proof}
Assume $G_T$ lies in a disk $D$ in $\hat{T}$.  By a sequence of disk exchanges, we may assume the interiors of the two faces of $G_S$ are disjoint from $D$.  Then $M=N(D \cup H_{1,2} \cup H_{2,1})$ is a genus 2 handlebody.  

The two faces of $G_S$ are each faces of order 3 Scharlemann cycles.  Attaching the $2$--handle neighborhoods of these two faces to $M$ will produce a connect sum of two lens spaces each of order 3.  This cannot be contained a lens space, however, since lens spaces are irreducible.
\end{proof}

\begin{figure}
\centering
\input{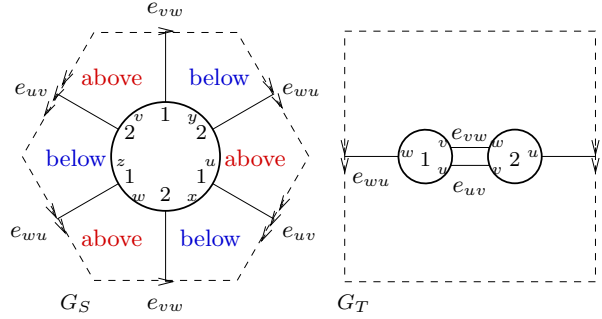}
\caption{$G_S$ of Figure~\ref{fig:feweredges}(b) and $G_T$ with labelings of edges and corners.}
\label{fig:markingGSandGT}
\end{figure}

\begin{lemma}\label{lem:somegenusoneknots}
Up to homeomorphism, each lens space $L(9k-3,3k-2)$, $k \in \Z$, contains exactly one order $3$ \opt--knot, $K^3_k$.  Every order $3$ \opt--knot is one of these.
\end{lemma}

\begin{proof}
By Lemma~\ref{lem:no3paralleledges} and Lemma~\ref{lem:no2paralleledges}, a non-torus knot of order $3$ with a properly embedded once-punctured torus in its exterior must have $G_S$ appear as in Figure~\ref{fig:feweredges}(b).
As shown, $G_S$ contains a single Scharlemann cycle of length $3$  $\sigma$ which bounds two faces $f_a$ and $f_b$ of colors amber and blue respectively.   

Lemma~\ref{lem:notinadisk} together with Lemma~\ref{lem:gt2.1} implies that $G_T$ lies in an essential annulus.  With labelings of the edges and corners of $f_a$ as in Figure~\ref{fig:meridian}, $G_S$ and $G_T$ may be assumed to be as shown in Figure~\ref{fig:markingGSandGT}.   Since handlebodies are irreducible, these labelings then dictate the positions of the disks $f_a$ and $f_b$ in the genus $2$ handlebodies $(V_a - H_{1,2}) \cup \bdry_v H_{1,2}$ and $(V_b - H_{2,1}) \cup \bdry_v H_{2,1}$ respectively.  
A priori the position of $\bdry f_a$ is determined up to Dehn twists along the core curve $h_a$ of $\bdry_v H_{1,2}$.  But since $h_a$ is dual to a compressing disk of the handlebody and $\bdry f_a$ intersects $h_a$ algebraically (and geometrically) $3$ times, there is at most one configuration for which $\bdry f_a$ bounds a disk.  Such a configuration is depicted in Figure~\ref{fig:length3schcycleinsolidtorus}.  The same argument holds for $f_b$.

\begin{figure}
\centering
\input{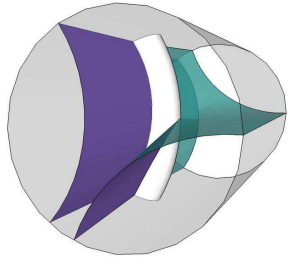}
\caption{Joining the front to the back with a $2\pi/3$ rotation forms the trigon face of a length $3$ Scharlemann cycle in a solid torus.}
\label{fig:length3schcycleinsolidtorus}
\end{figure}

Gluing $V_a$ and $V_b$ together along $\hat{T}$ so that the edges of $f_a$ and $f_b$ meet correctly on $G_T$ reconstitutes our lens space, our knot $K$, and once-punctured torus $S$.  Observe that since $S \cup K$ is disjoint from an essential annulus of $\hat{T}$, there is a $1$--parameter family of possible lens spaces formed by Dehn twists along this annulus of the gluing map.  To determine what lens spaces are thusly produced, we examine how the amber and blue meridians are situated on $\hat{T}$ with respect to one another.

With the given labelings, Lemma~\ref{lem:meridiansandS3cycles} implies that the amber and blue meridians, $m_a$ and $m_b$ respectively, lie on $\hat{T}$ with respect to $G_T$ as shown in Figure~\ref{fig:bothmeridiansGS} each up to Dehn twists along an essential simple closed curve $C$ on $\hat{T}$ that is disjoint from $G_T$.  Since both meridians are determined up to Dehn twists along $C$, we may fix $m_b$ and twist $m_a$.  Let $\tau_C(m_a)$ denote a single Dehn twist of $m_a$ along $C$ to the right.

\begin{figure}
\centering
\input{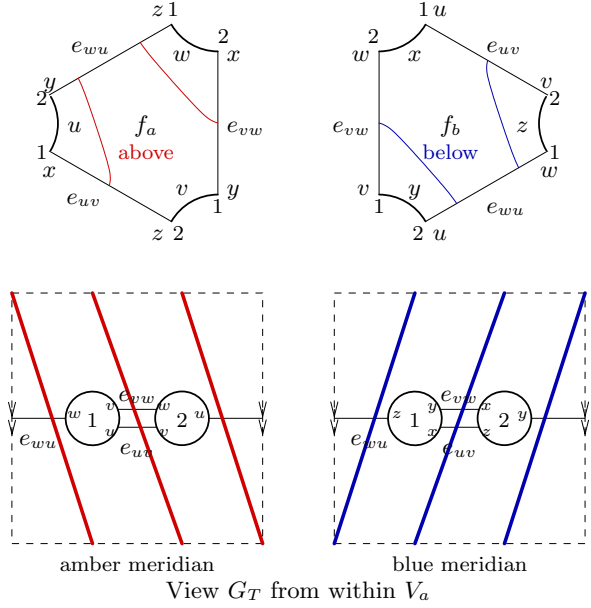}
\caption{The amber and blue meridional disks intersect the amber and blue faces, $f_a$ and $f_b$ respectively.  The amber and blue meridians then lie on $\hat{T}$ with respect to $G_T$ as shown.}
\label{fig:bothmeridiansGS}
\end{figure}

In Figure~\ref{fig:bothmeridiansGS}, orient $m_b$ upwards and rightwards.  Let $l$ be a vertical curve oriented upwards.  Then $\{[m_b],[l]\}$ forms a basis for $H_1(\hat{T})$.  Orient $m_a$ downwards and rightwards; orient $C$ rightwards.  Then $[m_a] = [m_b]- 6[l]$ and $[C] = [m_b]-3[l]$.  Since $\Delta(m_a,C) = 3$, $[\tau_C^k(m_a)] = [m_a]+3k[C] = (1+3k)[m_b] - (6+9k)[l]$.  Therefore the lens spaces that may be obtained from this construction with $G_T$ as shown are $L(-(6+9k),1+3k)$ for $k \in \Z$.  We make the reparametrization $k \mapsto k-1$ and apply homeomorphisms to obtain $L(9k-3,3k-2)$.  Then let $K^3_k$ be the knot $K$ in $L(9k-3,3k-2)$.
\end{proof}

\begin{remark} 
Since $\Delta(m_a, C)=3$ and $\Delta(m_b,C)=3$, the curve $C$ is a torus knot with exterior $W(-3, \cdot)$; see Lemma~\ref{lem:torusknots}.  The order $3$ \opt--knots $K^3_k$ in the lens spaces $L(9k-3,3k-2)$ are related by Dehn surgery on $C$.
\end{remark}

\begin{lemma}\label{lem:order3gridnumber}
The order $3$ \opt--knot $K^3_k$ is the $|3k-1|$th grid number $1$ knot in the lens space $L(9k-3,3k-2)$.
\end{lemma}

\begin{proof}
Write $K^3_k = K_a \cup K_b$.  Recall $K_a$ and $K_b$ are each unknotted arcs in $V_a$ and $V_b$ respectively.  By Lemma~\ref{lem:meridiansandS3cycles}  there are meridians $m_a$ and $m_b$ of $V_a$ and $V_b$ that bound meridional disks that are disjoint from each $K_a$ and $K_b$ and intersect $G_T$ in a prescribed manner.  Attaching $\bdry V_a$ to $\bdry V_b$ so that they agree on $G_T$ may be done so that $m_a$ and $m_b$ intersect minimally.  Since the vertices of $G_T$ mark where $K$ meets $\hat{T}$ and each vertex of $G_T$ is in its own parallelogram of the toroidal grid $m_a \cup m_b$, $K$ is grid number $1$.  Since $K$ has order $3$, $K$ must be the $|r/3|$th grid number $1$ knot, up to homeomorphism.
\end{proof}

\begin{lemma}\label{lem:order3surgdesc}
The order $3$ \opt--knot $K^3_k$ in the lens space $L(9k-3,3k-2)$ is the core of the $-3$--surgery in $W(-3,-3+1/k)$.
\end{lemma}

\begin{proof}
By Lemma~\ref{lem:surgeryrelationships} $W(-3,-3+1/k) = N(-2,-2+1/k,1)$.  Martelli and Petronio identify this as the lens space $L(9k-3,3k-2)$ in Table 11 of \cite{mp:dfotm3m}.  

Let $K'_k$ be the core of the $3$--surgery in $W(3,3-1/k)$.  Since the component of $W$ with the surgery coefficient $3$ bounds a once-punctured torus $S$ that is disjoint from the other component $C$ of $W$,  varying the surgery coefficient on $C$ preserves the surface $S$ and the meridian of the core of the $3$--surgery remains distance $3$ from $\bdry S$.  Hence for each integer $k$, $K'_k$ is a knot of order $3$ in the lens space $L(9k-3,3k-2)$ with the once-punctured torus $S$ properly embedded in its exterior.  Therefore by Lemma~\ref{lem:somegenusoneknots} the knot $K'_k$ must be the knot $K^3_k$.
\end{proof}

\begin{lemma}\label{lem:somegenusonefibered}
The exterior of $K^3_k$ in the lens space $L(9k-3,3k-2)$, $k \in \Z$, is fibered with fiber $S$ if an only if $k=\pm 1$.  Furthermore $S$ compresses if and only if $k=0$ in which case the knot exterior is a solid torus and $K^3_0$ is unknotted.
\end{lemma}

\begin{proof}
Since $K^3_k$ is the core of the $-3$--surgery in $W(-3,-3+1/k)$, its exterior is $W(-3+1/k, \cdot)$.  By Lemma~\ref{lem:whiteheadfiber} this exterior is fibered if and only if $-3+1/k \in \Z$, i.e.\ $k = \pm1$.  

If $k=0$ then $K^3_0$ is the core of $3$ surgery on the unknot in $S^3$ and hence its exterior is a solid torus and $S$ compresses.  If for some $k\neq0$ the surface $S$ compresses to a disk $S'$ then the neighborhood of $K^3_k \cup S$ would be a punctured order $3$ lens space in a lens space of order other than $3$; this cannot occur.
\end{proof}

\begin{lemma}\label{lem:order3fiber}
The knot $K^3_{+1}$ is a type $\{2,4\}$ torus knot in $L(6,1)$.  Its exterior is homeomorphic to $W(-2,\cdot)$. 
\end{lemma}
\begin{proof}
Write $K=K^3_{+1}$.
Following Figure~\ref{fig:bothmeridiansGS}, the meridians $m_a$ and $m_b$ of amber and blue meridional disks $D_a$ and $D_b$, respectively, that are disjoint from $K$ in $L(6,1)$ appear on $\hat{T}$ with respect to $G_T$ as in Figure~\ref{fig:whenk=-1}(i).  Then $K_a$ and $K_b$ may be isotoped rel--$\bdry$ onto $\hat{T}$ as in Figure~\ref{fig:whenk=-1}(ii).  Hence $K$ is a torus knot.  
\begin{figure}
\centering
\input{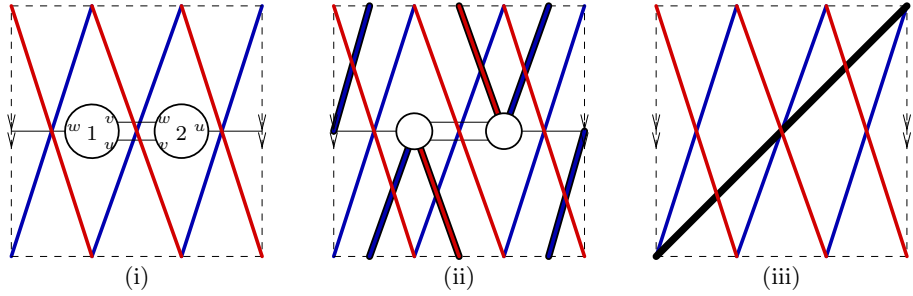}
\caption{An isotopy of $K^3_{+1}$ onto $\hat{T}$ in $L(6,1)$.}
\label{fig:whenk=-1}
\end{figure}
After ``straightening'' with an isotopy on $\hat{T}$ as in Figure~\ref{fig:whenk=-1}(iii), it is apparent that $\Delta(K,m_a)=4$ and $\Delta(K,m_b)=2$ and $K$ is a type $\{2,4\}$ torus knot.  Therefore, as in Lemma~\ref{lem:torusknots}, $E(K)$ is the torus knot exterior  $W(-2,\cdot)$.
\end{proof}

\begin{lemma}\label{lem:order3kplus1}
The knot $K^3_{-1}$ is a non-torus knot in $L(12,5)$.  Its exterior is the once-punctured torus bundle $W(-4, \cdot)$ with monodomy $\phi = \tau_x^4 \tau_y$.
\end{lemma}

\begin{proof}
Since the exterior of $K^3_{-1}$ is homeomorphic to $W(-4,\cdot)$, by Lemma~\ref{lem:whiteheadfiber} the exterior of $K^3_{-1}$ is fibered with monodromy $\phi \cong \tau_x^4 \tau_y$.
\end{proof}

\begin{remark}\label{rem:suturedmanifold}
Lemma~\ref{lem:somegenusoneknots} and Lemma~\ref{lem:somegenusonefibered} may be proven quite explicitly and graphically using the sutured manifold decomposition (see \cite{gabai:fatto3m}) of $(E(K), \bdry E(K))$, $K=K^3_k$.   First chop $V_a \cap E(K)$ and $V_b \cap E(K)$ along the amber and blue trigon faces of $S$ and then reattaching the resulting pieces along the remnants of $\hat{T}$.  This leaves one with a sutured manifold that is a genus $2$ handlebody $M$ with an annular suture $A$ that divides $\bdry M$ into two once-punctured tori.  In particular $(M,A)$ is the result of the sutured manifold decomposition of $(E(K), \bdry E(K))$ along $S$.  

When $(M,A)$ admits a complete product disk decomposition, $S$ is the fiber of a fibration of $E(K)$.  The monodromy of this fibration may be inferred from this product disk decomposition.  One may thusly recover Lemma~\ref{lem:order3fiber} and Lemma~\ref{lem:order3kplus1}.

We may now continue in this manner to explicitly see the surgery description of $K^3_k$ of Lemma~\ref{lem:order3surgdesc}.   One of the remnants of $\hat{T}$ is an annulus whose core is the curve $C$, and therefore lies in the interior of $M$.   In fact one can see that $C$ transversally intersects a compressing disk of $M$ once and is parallel to a curve on $\bdry M$ that crosses $A$ twice minimally.   Moreover, $(M-N(C), A)$ is homeomorphic to $(M'-N(C'),A')$ where $M'$ is the genus $2$ handlebody closed neighborhood of a theta graph $C' \cup a'$ formed from attaching the endpoints of an arc $a'$ to distinct points on a circle $C'$ and $A'$ is an essential annulus on $\bdry M'$ whose core bounds a disk in $M$ that intersects $C' \cup a'$ only in $a'$.  One may now identify $(M'-N(C'), A')$ with the exterior of $W$ chopped along a properly embedded once-punctured torus.
\end{remark}

\section{$s=2$}\label{sec:s=2}
This section is devoted to proving the following theorem.
\begin{thm} \label{thm:seq2}
If $s=2$ then, for each integer $k \neq 0$, $K$ is a knot in $W(-2,-4+1/k) \cong L(8k-2,4k+1)$ arising as the core of the $2$ surgery.  Equivalently, up to homeomorphism, $K$ is the $|4k-1|$th grid number $1$ knot in the lens space $L(8k-2,4k+1)$.  Furthermore $E(K)$ is fibered if and only if $k = \pm 1$:
\begin{itemize}
\item If $k=+1$ then $K$ is a torus knot in $L(6,1)$ with monodromy $\phi \cong \tau_x^3 \tau_y$.
\item If $k=-1$ then $K$ is a non-torus knot in $L(10,3)$ with monodromy $\phi \cong \tau_x^5 \tau_y$.
\end{itemize}
\end{thm}

\begin{proof}
We give an overview of the proof.  The first two subsections culminate in Theorem~\ref{thm:s=2t=2} showing that if $s=2$ then $t=2$ and $G_S$ appears as in Figure~\ref{fig:G_Ssis2}.  Then Lemma~\ref{lem:somegenusoneknotsorder2} describes the lens spaces which contain knots with such surfaces $S$.  Lemma~\ref{lem:order2surgdesc} gives the surgery description of this family of knots.   Lemma~\ref{lem:somegenusonefiberedorder2} shows that $E(K)$ is fibered if and only if $k = \pm1$ and that $S$ is an incompressible if and only if $k \neq 0$.  Lemma~\ref{lem:order2kminus1} shows that if $k=+1$ then $K$ is a torus knot. Lemma~\ref{lem:order2kplus1} shows that if $k=-1$ then $K$ has fibered exterior with monodromy $\phi \cong \tau_x^5 \tau^y$.
\end{proof}

If $s=2$ then Theorem~1.1 of \cite{baker:sgkilshsbn} give no apriori bounds on $t$.  Nevertheless, we will explicitly conclude that here too, $t=2$.

When $s=2$, the boundary of $S$ runs twice longitudinally along $\bdry N(K)$.
Within $N(K)$ there is a Mobius band $P$ with the same boundary as $S$.  Furthermore we may take $K$ to be the core curve of $P$.  Then $S \cup P$ is homeomorphic to $\hat{S} \# \hat{P}$ (where $\hat{P} \cong \R P^2$ thought of as abstractly capping off $P$ with a  disk).  Indeed for every closed non-orientable surface $R$ of Euler characteristic $-1$ embedded in the lens space $X$ and orientation reversing simple closed curve $J \subset R$ such that $R - N(J)$ is a once-punctured torus incompressible in $E(J)$, we have a knot $J$ with an incompressible once-punctured torus Seifert surface $R-N(J)$.

\begin{lemma}\label{lem:singlecurve}
There is a unique isotopy class of simple closed curves $J$ in $R \cong \#_3 \R P^2$ such that $J$ is orientation reversing and $R - N(J)$ is a once-punctured torus.
\end{lemma}

\begin{proof}
We leave this as a fun exercise.
\end{proof}

\subsection{When $s=2$ and $r=2$}

\begin{lemma}\label{lem:nogenus1inL(2,1)}
An \opt--knot in $L(2,1)$ of order $2$ is unknotted.
\end{lemma}

\begin{proof}
Let $K$ be a knot of order $2$ in $L(2,1)$ with a once-punctured torus $S$ properly embedded in its exterior.    Extend $S$ radially through $N(K)$ to form the surface $R \cong \#_3 \R P^2$ embedded in $L(2,1)$.  

Identify $V_a$ with $S^1 \times D^2$ and $V_b$ with $D^2 \times S^1$.  We proceed in the vein of Bredon-Wood \cite{bw:nosio3m}.  Isotop $R$ to intersect $S^1 \times \{0\}$ transversally.  $R$ must intersect $S^1 \times \{0\}$ an odd number of times, $d$.  Re-envision $\hat{T}$ as the boundary of a small tubular neighborhood of this circle and adjust $V_a$ and $V_b$ accordingly.  Then $R$ intersects $V_a$ in $d$ meridional disks.  The boundary of $R \cap V_b$ is then $d$ parallel curves  of slope $1/2$ on $\bdry V_b$.  We will now isotop $R$ to make $d=1$.

Isotop $R$ to meet the meridional disk $D = D^2 \times \{1\}$ of $V_b$ transversally. Then $R \cap D$ consists of $d$ properly embedded arcs and a collection of simple closed curves.  Disregarding simple closed curves, if $d >1$ then an outermost arc of $R \cap D$ on $D$ cuts off a unique disk $E_b$ that is disjoint from the other arcs.  The two curves on $\bdry V_b$ that $E_b$ meets each bound meridional disks of $V_a$ so that there is a product interval $I \times D^2$ between them that is disjoint from the rest of $R$ and $E_b \cap I \times D^2$ is a spanning arc of the annulus $I \times \bdry D^2$.  Let $E_a$ be a radial extension of this spanning arc to the arc $I \times \{0\}$. Use $E = E_a \cup E_b$ to isotop a neighborhood of $E$ in $R$ (including any simple closed curves) to the ``other side'' of  $S^1 \times \{0\} \subset V_a$ by pushing from $R \cap \bdry E$ through $E$ to just past $\bdry E - R$.  This reduces the number of times $R$ intersects $S^1 \times \{0\}$ by 2.  Rechoose $\hat{T}$ to be the boundary of (a yet smaller) tubular neighborhood of $S^1 \times \{0\}$ and adjust $V_a$ and $V_b$ accordingly.  Then $R$ intersects $V_a$ in $d-2$ meridional disks.  Repeat this process until $d=1$.

Now that $d=1$, $R \cap V_b$ may be seen to be a properly embedded M\"obius band $P$ with a $1$--handle attached.  We may isotop the feet of this $1$--handle into a small disk $\Delta$ in $P$ such that both feet are attached to the same side of $\Delta$.  Let $J$ be a core curve of $P$ that is disjoint from $\Delta$.  Then $J$ is an orientation reversing curve in $R$ and $R - N(J)$ is a once-punctured torus.  By Lemma~\ref{lem:singlecurve} $J$ is the unique simple closed curve on $R$ up to isotopy with these two properties.  Therefore $J$ is isotopic on $R$ to $K$.  One now observes that $J$, and hence $K$, is isotopic to the core of $V_b$.  Therefore the exterior of $K$ is a solid torus and $K$ is unknotted.  (Furthermore the co-core of the $1$--handle is a compressing disk for $R$ disjoint from $J$ thereby giving a compression of $S$.)
\end{proof}

\subsection{When $s=2$ and $r>2$...}
In light of Lemma~\ref{lem:nogenus1inL(2,1)} we assume $r >2$.  Since $r \neq 2$,  Lemmas~\ref{lem:2.2ish} and \ref{lem:gt2.3} imply that the edges of a Scharlemann cycle of length two of $G_S$ must lie in an essential annulus on $G_T$.
  
\begin{figure}
\centering
\input{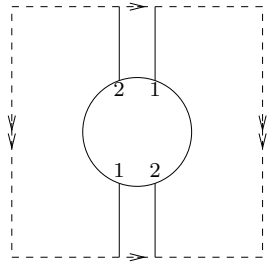}
\caption{The only viable configuration for $G_S$ when $s=2$. (Simple closed curves not included.)}
\label{fig:G_Ssis2}
\end{figure}

\begin{thm}\label{thm:s=2t=2}
If $s=2$ then $t=2$ and $G_S$ appears as in Figure~\ref{fig:G_Ssis2}.
\end{thm}

\begin{proof}
Assume $t > 2$.
\begin{figure}
\centering
\input{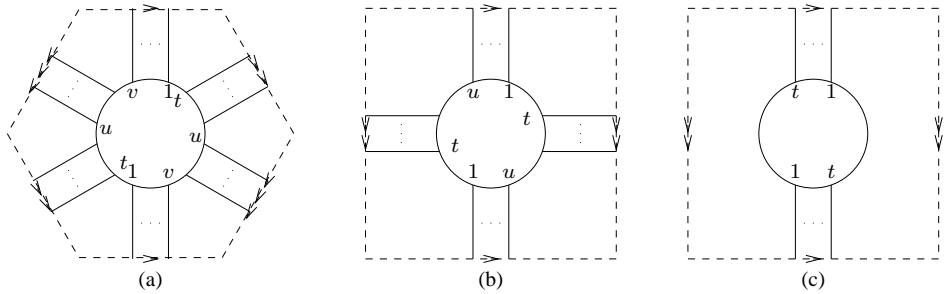}
\caption{When $s=2$ and $G_S$ has three, two, or one edge classes.}
\label{fig:G_Sfors=2}
\end{figure}
Around the vertex of $G_S$ the labels $\{1, \dots, t\}$ each appear twice.  By relabeling if necessary, $G_S$ appears as in Figure~\ref{fig:G_Sfors=2} (a), (b), or (c) according to whether $G_S$ has three, two, or one edge classes respectively.  Below, Lemma~\ref{lem:3edgeclassess=2} and Lemma~\ref{lem:2edgeclassess=2} show that $G_S$ must have just one edge class and Lemma~\ref{lem:1edgeclassess=2} shows that $t$ cannot be greater than $2$.  The configuration of $G_S$ is then forced.
\end{proof}

\begin{lemma}\label{lem:3edgeclassess=2}
If $s=2$ then $G_S$ cannot have three edge classes.
\end{lemma}
\begin{proof}
Assume $G_S$ has three edge classes. 
\begin{figure}
\centering
\input{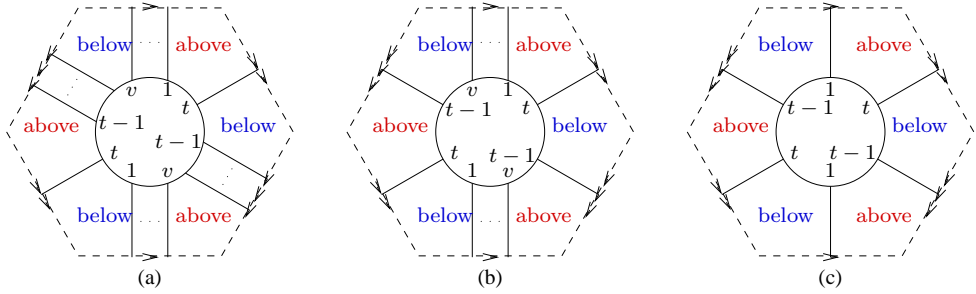}
\caption{When $s=2$, $G_S$ has three edge classes, and some edge class has just one edge.}
\label{fig:G_Ss=23edgeclasses}
\end{figure} 
Each edge class of $G_S$ must have at least two edges.  If not then $G_S$ appears up to relabeling as in Figure~\ref{fig:G_Ss=23edgeclasses}(a), (b), or (c) according to whether just one, two, or three edge classes have just one edge respectively.  In each case there are two distinct faces, $f_1$ and $f_2$, each with a corner on the blue $1$--handle $H_{t, 1}$.  Since any two faces meeting along an edge have opposite colors, any choice of coloring of a trigon of $G_S$ in turn forces $f_1$ and $f_2$ to have opposite colors.  Yet because $f_1$ and $f_2$ both have a corner on $H_{t,1}$, they must both be blue; a contradiction.

Now since each edge class has at least two edges then, as evident from Figure~\ref{fig:G_Sfors=2}(a), $G_S$ has three Scharlemann cycles of length two with mutually disjoint label pairs.  By Lemma~\ref{lem:gtS2cycles} this cannot occur.  
\end{proof}

\begin{lemma}\label{lem:2edgeclassess=2}
If $s=2$ then $G_S$ cannot have two edge classes.
\end{lemma}
\begin{proof}
Assume $G_S$ has two edge classes as in Figure~\ref{fig:G_Sfors=2}(b).  
\begin{figure}
\centering
\input{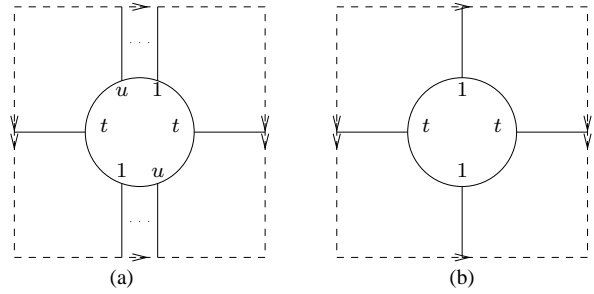}
\caption{When $s=2$, $G_S$ has two edge classes, and some edge class has just one edge.}
\label{fig:G_Ss=22edgeclasses}
\end{figure} 
Each edge class of $G_S$ must have at least two edges.  If not then $G_S$ appears up to relabeling as in Figure~\ref{fig:G_Ss=22edgeclasses}(a) or (b) according to whether just one or two edge classes have just one edge respectively.  Both of these situations fail to satisfy the Parity Rule.
 
Then as evident from Figure~\ref{fig:G_Sfors=2}(b), $G_S$ has two Scharlemann cycles of length two, $\sigma$ and $\tau$, with disjoint label pairs.  By Lemma~\ref{lem:gtS2cycles} the faces of $\sigma$ and $\tau$ must have different colors.  

There is a tetragon face $\rho$ of $G_S$.  Let $f$ be the face bounded by $\rho$.  We claim that $\Int f \cap \hat{T}$ does not contain any curves that are essential on $\hat{T}$.

If not then $\Int f \cap \hat{T}$ contains a curve that is essential on $\hat{T}$.  After any necessary disk exchanges, there is an innermost disk $D$ on $f$ bounded by an essential curve of $\Int f \cap \hat{T}$.  One of the Scharlemann cycles, say $\sigma$, has the same color as $D$.  Joining the corners of $\sigma$ along $K$ forms a M\"obius band whose boundary is disjoint and parallel to $\bdry D$ on $\hat{T}$.  Then with the annulus of parallelism, the M\"obius band and $D$ form a projective plane in a solid torus.  This cannot occur.  

Given the claim, by disk exchanges we may arrange that $\Int f \cap \hat{T} = \emptyset$.   Opposite corners of $f$ have the same two labels, either $\{t, 1\}$ or $\{u, u+1\}$.  Extend the corners of $f$ radially to the cores of $H_{t,1}$ and $H_{u,u+1}$ to form the twice punctured projective plane $f'$.
Since the opposite edges of $f$ each bound extended Scharlemann cycles in $G_S$, the boundary components of $f'$ lie in essential annuli and are therefore parallel.  Let $A$ be this annulus of parallelism.  Then $f' \cup A \cong \R P^3 \# T^2$ is a closed non-orientable surface in a solid torus.  This cannot occur.
\end{proof}

\begin{lemma}\label{lem:1edgeclassess=2}
If $s=2$ and $t>2$ then $G_S$ cannot have one edge class.
\end{lemma}

\begin{proof}
Assume $G_S$ has one edge class.  Then as evident from Figure~\ref{fig:G_Sfors=2}(c), $G_S$ has an extended Scharlemann cycle of length two, $\sigma$, that contains all of the labels.  (Note that this edge class cannot have just one edge.)  Let $g$ be the bigon face bounded by the Scharlemann cycle in $\sigma$.  

Since $\sigma$ pairs edges of $G_T$ so that the pairs lie in essential annuli on $\hat{T}$, no component of $G_T$ lies in a disk on $\hat{T}$.  Therefore we may perform disk exchanges to eliminate simple closed curves of $S \cap T$ that are trivial on $\hat{T}$.  Then the bigons of the extended Scharlemann cycle give rise to a Mobius band and a collection of annuli by extending their corners radially to the core of $H_{1,t} = H_{1,2} \cup \dots \cup H_{t-1,t}$.  The bigon $g$ gives a central Mobius band $A_0$.  Stepping uniformly away from $g$, the $i$th pair of bigons form an annulus $A_i$ which shares a boundary component for with $A_{i-1}$ for $i=1, \dots, t/2$.  Hence we form the long Mobius band $A = A_0 \cup A_1 \cup \dots A_{t/2}$.  Since $K_{t,1} = H_{t,1} \cap K$ is the only arc of $K-\hat{T}$ not contained in $A$, we may use $A$ to thin $K$ unless $t=2$.  See $\S4$ of \cite{baker:sgkilshsbn}. 
\end{proof}

\subsection{When $G_S$ is as in Figure~\ref{fig:G_Ssis2}...}

\begin{lemma}\label{lem:isotopsis2}
If $s=2$ then $S \cap T$ contains no simple closed curves.
\end{lemma}

\begin{proof}
By Theorem~\ref{thm:s=2t=2}, $G_S$ appears as in Figure~\ref{fig:G_Ssis2}.  Let $g$ be the disk face bounded by the Scharlemann cycle.  Let $f$ be the annulus face.

Since $p \neq 2$ (Lemma~\ref{lem:nogenus1inL(2,1)}) the two edges of $G_T$ lie in an essential annulus (Lemmas~\ref{lem:2.2ish} and \ref{lem:gt2.3}).  By disk exchanges we may assume any simple closed curve of $S \cap T$ is essential on $\hat{T}$.  Since lens spaces are irreducible and $K$ is not contained in a ball, these disk exchanges may be realized by isotopies of $S$ with support in the exterior of $K$.

If $S \cap T$ contains a simple closed curve that bounds a disk on $S$, then assume $D$ is an innermost such disk.  By Lemma~\ref{lem:gt2.3} this disk must lie in the annulus face $f$ of $G_S$.  Form a Mobius band $A_0$ by extending the corners of $g$ radially to the core of $H_{1,2}$.  Then since $\bdry A_0$ and $\bdry D$ are disjoint essential curves on $\hat{T}$, they are parallel.  Let $B$ be an annulus of parallelism between these curves on $\hat{T}$.  Then $N(B \cup A_0 \cup D)$ is a punctured $L(2,1)$.  Since $p \neq 2$, this is a contradiction.  Hence any simple closed curve of $S \cap T$ is essential on $S$.  Such curves are all contained in the annulus face $f$ of $G_S$, and there must be an even number of them.  

We isotop $S$ to reduce the number of these curves.
Include all the simple closed curves of $S \cap T$ in $G_S$ and color the faces so that $g$ (and the corresponding M\"obius band $A_0$) is colored blue.  Then every blue annulus is disjoint from $A_0$ and may be regarded as the horizontal boundary of a twisted $I$ bundle over $A_0$ in $V_b$.  Moreover, each blue annulus is boundary-parallel into an annulus of $\hat{T}$ through the solid torus complement of its corresponding twisted $I$ bundle in $V_b$.  Isotop an outermost blue annulus onto $\hat{T}$ and then push it slightly into $V_a$.  This reduces the number of simple closed curves of $S \cap T$ by two.  

Therefore if $S \cap T$ contains a simple closed curve then $|S \cap T|$ is not minimized.  This is a contradiction to our minimality assumption on $|S \cap T|$.  Hence $S \cap T$ contains no simple closed curves.
\end{proof}

\begin{figure}
\centering
\input{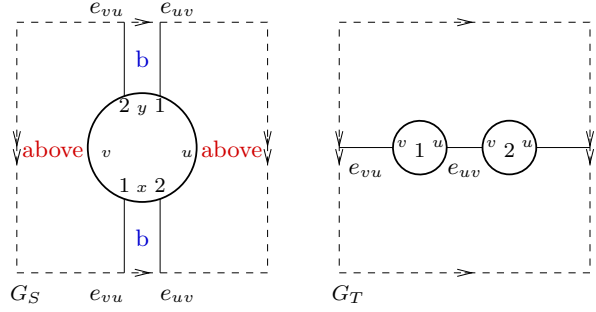}
\caption{$G_S$ of Figure~\ref{fig:G_Ssis2} and $G_T$ with labelings of edges and corners.}
\label{fig:markingGSandGTorder2}
\end{figure}

\begin{lemma}\label{lem:somegenusoneknotsorder2}
Up to homeomorphism, each lens space $L(8k-2,4k+1)$, $k\in \Z$, contains exactly one order $2$ \opt--knot $K^2_k$.  Every order $2$ \opt--knot is one of these.
\end{lemma}

\begin{proof}
By Theorem~\ref{thm:s=2t=2}, a non-torus \opt-knot of order $2$ must have $G_S$ appear as in Figure~\ref{fig:G_Ssis2}.  As shown, $G_S$ contains a single Scharlemann cycle of length $2$ which bounds an annulus face $f_a$ and a bigon face $f_b$ of colors amber and blue respectively.

We label the edges and corners of $f_a$ and $f_b$ so that $G_S$ and $G_T$ may be assumed to be as shown in Figure~\ref{fig:markingGSandGTorder2}.  

By extending the corners of $f_a$ radially to the core of $H_{1,2}$ we obtain a once-punctured Klein bottle $\bar{f}_a$ properly embedded in $V_a$ that does not compress to a M\"obius band.  Therefore $\Delta(\bdry \bar{f}_a, m_a)=4$, and up to homeomorphism $f_a$ appears as in Figure~\ref{fig:length2schcycleinsolidtorus}(a).  Similarly $f_b$ may be extended to a M\"obius band $\bar{f}_b$ properly embedded in $V_b$.  Then $\Delta(\bdry \bar{f}_b,m_b)=2$, and $f_b$ appears as in Figure~\ref{fig:length2schcycleinsolidtorus}(b).  Observe that $K_a$ and $K_b$ are disjoint from meridional disks of $V_a$ and $V_b$ bounded by $m_a$ and $m_b$ respectively.

\begin{figure}
\centering
\input{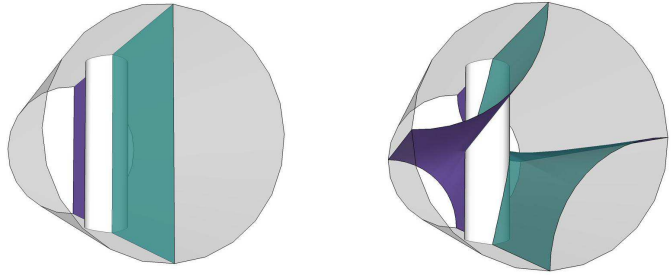}
\caption{(a) Joining the front to the back with a $\pi/2$ rotation forms the annulus $f_a$ in the solid torus $V_a$.  (b) Joining the front to the back with a $\pi$ rotation forms the bigon face $f_b$ of a length $2$ Scharlemann cycle in a solid torus.}
\label{fig:length2schcycleinsolidtorus}
\end{figure}

Gluing $V_a$ and $V_b$ together along $\hat{T}$ so that the edges of $f_a$ and $f_b$ meet correctly on $G_T$ reconstitutes our lens space, our knot $K$, and once-punctured torus $S$.  Observe that since $S \cup K$ is disjoint from an essential annulus of $\hat{T}$, there is a $1$--parameter family of possible lens spaces formed by Dehn twists along this annulus of the gluing map.  To determine what lens spaces are thusly be produced, we examine how the amber and blue meridians are situated on $\hat{T}$ with respect to one another.

Observe that $\bdry \bar{f}_a$ and $\bdry \bar{f}_b$ are $G_T$ with the vertices contracted to a point on $\hat{T}$.   Then since $\Delta(\bdry \bar{f}_a,m_a)=4$ and $\Delta(\bdry \bar{f}_b,m_b)=2$, up to homeomorphism, $m_a$ and $m_b$ lie on $\hat{T}$ with respect to $G_T$ as shown in Figure~\ref{fig:bothmeridiansGSorder2} each up to Dehn twists along an essential simple closed curve $C$ on $\hat{T}$ that is disjoint from $G_T$.  Since both meridians are determined up to Dehn twists along $C$, we may fix $m_b$ and twist $m_a$.  Let $\tau_C(m_a)$ denote a single Dehn twist of $m_a$ along $C$ to the right.

\begin{figure}
\centering
\input{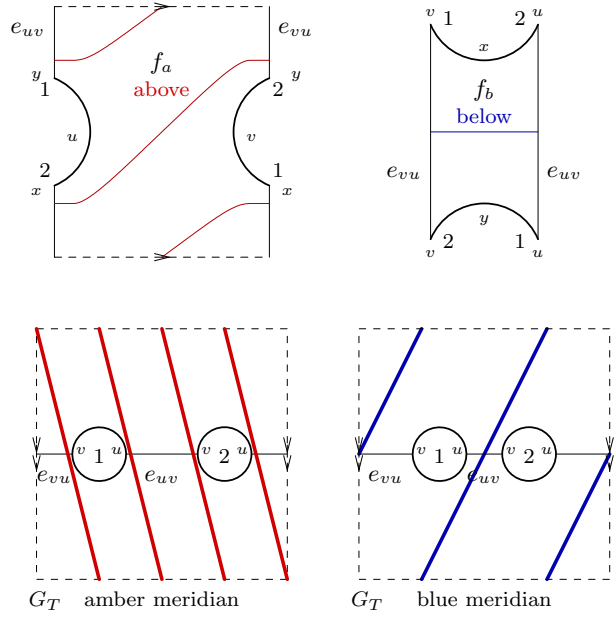}
\caption{The amber and blue meridional disks intersect the amber and blue faces, $f_a$ and $f_b$ respectively.  The amber and blue meridians then lie on $\hat{T}$ with respect to $G_T$ as shown.}
\label{fig:bothmeridiansGSorder2}
\end{figure}

In Figure~\ref{fig:bothmeridiansGS}, orient $m_b$ upwards and rightwards.  Let $l$ be a vertical curve oriented upwards.  Then $\{[m_b],[l]\}$ forms a basis for $H_1(\hat{T})$.  Orient $m_a$ downwards and rightwards; orient $c$ rightwards.  Then $[m_a] = [m_b]- 6[l]$ and $[C] = [m_b]-2[l]$.  Since $\Delta(m_a,C)=4$, $[\tau_C^k(m_a)]=[m_a]+4k[C] = (1+4k)[m_b]-(6+8k)[l]$.  Therefore the lens spaces that may be obtained from this construction with $G_T$ as show are $L(-(6+8k), 1+4k)$.  We make the reparametrization $k \mapsto k-1$ and apply homeomorphisms to obtain $L(8k-2,4k+1)$.  Let $K^2_k$ be the knot $K$ in $L(8k-2,4k+1)$.
\end{proof}

\begin{remark} 
Since $\Delta(m_a, C)=4$ and $\Delta(m_b,C)=2$, the curve $C$ is a torus knot with exterior $W(-2, \cdot)$; see Lemma~\ref{lem:torusknots}.  The order $2$ \opt--knots $K^2_k$ in the lens spaces $L(8k-2,4k+1)$ are related by Dehn surgery on $C$.
\end{remark}

\begin{lemma}\label{lem:order2gridnumber}
The order $2$ \opt--knot $K^2_k$ is the $|4k-1|$th grid number $1$ knot in the lens space $L(8k-2,4k+1)$.
\end{lemma}

\begin{proof}
Write $K^2_k = K_a \cup K_b$.  Recall $K_a$ and $K_b$ are each unknotted arcs in $V_a$ and $V_b$ respectively.  As noted in the proof of Lemma~\ref{lem:somegenusoneknotsorder2}, there are meridians $m_a$ and $m_b$ of $V_a$ and $V_b$ that bound meridional disks that are disjoint from each $K_a$ and $K_b$ and intersect $G_T$ in a prescribed manner.  Attaching $\bdry V_a$ to $\bdry V_b$ so that they agree on $G_T$ may be done so that $m_a$ and $m_b$ intersect minimally.  Since the vertices of $G_T$ mark where $K$ meets $\hat{T}$ and each vertex of $G_T$ is in its own parallelogram of the toroidal grid $m_a \cup m_b$, $K$ is grid number $1$.  Since $K$ has order $2$, $K$ must be the $|r/2|$th grid number $1$ knot.
\end{proof}

\begin{lemma}\label{lem:order2surgdesc}
The knot $K^2_k$ in the lens space $L(8k-2,4k+1)$ is the core of the $-2$--surgery in $W(-2,-4+1/k)$.
\end{lemma}

\begin{proof}
This proof follows exactly the same as that of Lemma~\ref{lem:order3surgdesc}.  
By Lemma~\ref{lem:surgeryrelationships} $W(2,4-1/k) = N(-1,-3+1/k,1)$.  Martelli and Petronio identify this as the lens space $L(8k-2,4k+1)$ in Table 11 of  \cite{mp:dfotm3m}.  

Let $K'_k$ be the core of the $-2$--surgery in $W(-2,-4+1/k)$.  Since the component of $W$ with the surgery coefficient $-2$ bounds a once-punctured torus $S$ that is disjoint from the other component $C$ of $W$,  varying the surgery coefficient on $C$ preserves the surface $S$ and the meridian of the core of the $2$--surgery remains distance $2$ from $\bdry S$.  Hence for each nonzero integer $k$, $K'_k$ is an \opt--knot of order $2$ in the lens space $L(8k-2,4k+1)$ with the once-punctured torus $S$ properly embedded in its exterior.  
Therefore by Lemma~\ref{lem:somegenusoneknotsorder2} the knot $K'_k$ must be the knot $K^2_k$.
\end{proof}

\begin{lemma}\label{lem:somegenusonefiberedorder2}
The exterior of $K^2_k$ in the lens space $L(8k-2,4k+1)$, $k \in \Z$, is fibered with fiber $S$ if an only if $k=\pm1$.  Furthermore $S$ compresses if and only if $k=0$ in which case the knot exterior is a solid torus and $K^2_0$ is unknotted.
\end{lemma}

\begin{proof}
This proof follows similarly to that of Lemma~\ref{lem:somegenusonefibered}.  Since $K^2_k$ is the core of the $-2$--surgery in $W(-2,-4+1/k)$, its exterior is $W(-4+1/k, \cdot)$. By Lemma~\ref{lem:whiteheadfiber} this exterior is fibered if and only if $-4+1/k \in \Z$, i.e.\ $k=\pm1$.  

If $k=0$ then $K^2_0$ is the core of $-2$--surgery on the unknot in $S^3$.  Hence $S$ compresses and its exterior is a solid torus.  If for some $k \neq 0$ the surface $S$ compresses to a disk $S'$ then the neighborhood of $K^2_k \cup S'$ would be a punctured order $2$ lens space in a lens space of order other than $2$; this cannot occur.
\end{proof}

\begin{lemma}\label{lem:order2kminus1}
The knot $K^2_{+1}$ is a type $\{3,3\}$ torus knot in $L(6,1)$.  Its exterior is homeomorphic to  $W(-3, \cdot)$.
\end{lemma}

\begin{proof}
Write $K=K^2_{+1}$.
Following Figure~\ref{fig:bothmeridiansGSorder2}
the meridians $m_a$ and $m_b$ of amber and blue meridional disks $D_a$ and $D_b$ respectively that are disjoint from $K$ in $L(6,1)$ appear on $\hat{T}$ with respect to $G_T$ as in Figure~\ref{fig:whenk=-1ord2}(i).  Then $K_a$ and $K_b$ may be isotoped rel--$\bdry$ onto $\hat{T}$ as in Figure~\ref{fig:whenk=-1ord2}(ii).  Hence $K$ is a torus knot.  
\begin{figure}
\centering
\input{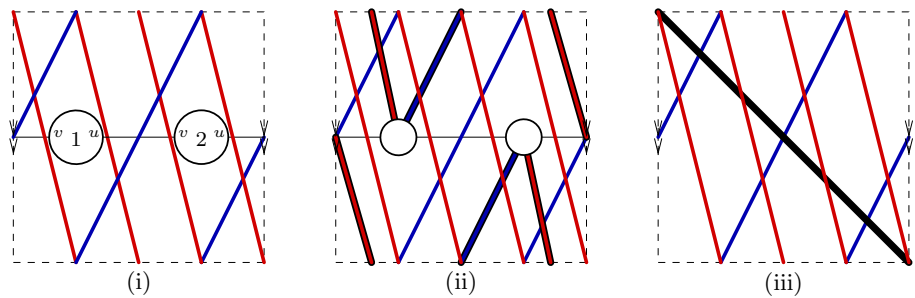}
\caption{An isotopy of $K^2_{+1}$ onto $\hat{T}$ in $L(6,1)$.}
\label{fig:whenk=-1ord2}
\end{figure}
After ``straightening'' with an isotopy on $\hat{T}$ as in Figure~\ref{fig:whenk=-1ord2}(iii), it is apparent that $\Delta(K,m_a)=3$ and $\Delta(K,m_b)=3$ and $K$ is a type $\{3,3\}$ torus knot.  Hence, as in Lemma~\ref{lem:torusknots}, $E(K)$ is the torus knot exterior $W(-3,\cdot)$.
\end{proof}

\begin{lemma}\label{lem:order2kplus1}
The knot $K^2_{-1}$ is a knot in $L(10,3)$ and its exterior is a once-punctured torus bundle with monodomy $\phi = \tau_x^5 \tau_y$.
\end{lemma}

\begin{proof}
Since the exterior of $K^2_{-1}$ is homeomorphic to $W(5,\cdot)$, by Lemma~\ref{lem:whiteheadfiber} the exterior of $K^2_{-1}$ is fibered with monodromy $\phi \cong \tau_x^5 \tau_y$.
\end{proof}

\begin{remark}
As in Remark~\ref{rem:suturedmanifold} we may here too use sutured manifolds to obtain explicit proofs of Lemma~\ref{lem:somegenusoneknotsorder2}, Lemma~\ref{lem:order2surgdesc}, Lemma~\ref{lem:somegenusonefiberedorder2}, Lemma~\ref{lem:order2kminus1}, and Lemma~\ref{lem:order2kplus1}.
\end{remark}


\begin{thebibliography}{HMW92}

\bibitem[Bak]{baker:cgofkils}
Kenneth~L. Baker.
\newblock {Counting genus one fibered knots in lens spaces},
  arXiv:math.GT/0510391.

\bibitem[Bak06]{baker:sgkilshsbn}
Kenneth~L. Baker.
\newblock Small genus knots in lens spaces have small bridge number.
\newblock {\em Algebr. Geom. Topol.}, 6:1519--1621 (electronic), 2006.

\bibitem[Bal]{baldwin:anogofkils}
John~A. Baldwin.
\newblock {A note on genus one fibered knots in lens spaces},
  arXiv:math.GT/0607370.

\bibitem[Ber]{berge:skwsyls}
John Berge.
\newblock Some knots with surgeries yielding lens spaces.
\newblock Unpublished manuscript.

\bibitem[BJK]{bjk:tnogofk}
Kenneth~L. Baker, Jesse~E. Johnson, and Elizabeth~A. Klodginski.
\newblock {Tunnel number one, genus one fibered knots}, arXiv:math.GT/0606379.

\bibitem[BO83]{bo:sdhdel}
Francis Bonahon and Jean-Pierre Otal.
\newblock Scindements de {H}eegaard des espaces lenticulaires.
\newblock {\em Ann. Sci. \'Ecole Norm. Sup. (4)}, 16(3):451--466 (1984), 1983.

\bibitem[BW69]{bw:nosio3m}
Glen~E. Bredon and John~W. Wood.
\newblock Non-orientable surfaces in orientable {$3$}-manifolds.
\newblock {\em Invent. Math.}, 7:83--110, 1969.

\bibitem[Bro87]{brown:tvatbnsi}
K.~S. Brown.
\newblock Trees, valuations, and the Bieri-Neumann-Strebel invariant.
\newblock {\em Invent. Math.}, 90(3):479Ð504, 1987.



\bibitem[BZ67]{bz:nkuf}
Gerhard Burde and Heiner Zieschang.
\newblock Neuwirthsche {K}noten und {F}l\"achenabbildungen.
\newblock {\em Abh. Math. Sem. Univ. Hamburg}, 31:239--246, 1967.

\bibitem[CB88]{cassonbleiler}
Andrew~J. Casson and Steven~A. Bleiler.
\newblock {\em Automorphisms of surfaces after {N}ielsen and {T}hurston},
  volume~9 of {\em London Mathematical Society Student Texts}.
\newblock Cambridge University Press, Cambridge, 1988.

\bibitem[CGLS85]{cgls:dsok}
Marc Culler, Cameron~McA. Gordon, John Luecke, and Peter~B. Shalen.
\newblock Dehn surgery on knots.
\newblock {\em Bull. Amer. Math. Soc. (N.S.)}, 13(1):43--45, 1985.

\bibitem[Cro98]{cromwell:apokal}
Peter~R. Cromwell.
\newblock Arc presentations of knots and links.
\newblock In {\em Knot theory (Warsaw, 1995)}, volume~42 of {\em Banach Center
  Publ.}, pages 57--64. Polish Acad. Sci., Warsaw, 1998.

\bibitem[Dyn06]{dynnikov:apolms}
I.~A. Dynnikov.
\newblock Arc-presentations of links: monotonic simplification.
\newblock {\em Fund. Math.}, 190:29--76, 2006.

\bibitem[Gab83]{gabai:fatto3m}
David Gabai.
\newblock Foliations and the topology of {$3$}-manifolds.
\newblock {\em J. Differential Geom.}, 18(3):445--503, 1983.

\bibitem[Gab87]{gabai:fatto3mIII}
David Gabai.
\newblock Foliations and the topology of {$3$}-manifolds. {III}.
\newblock {\em J. Differential Geom.}, 26(3):479--536, 1987.

\bibitem[GAn70]{ga:dcok}
Francisco Gonz\'alez-Acu\~na.
\newblock Dehn's construction on knots.
\newblock {\em Bol. Soc. Mat Mexicana}, pages 58--79, 1970.

\bibitem[GL95]{gl:dsokcetI}
C.~McA. Gordon and J.~Luecke.
\newblock Dehn surgeries on knots creating essential tori. {I}.
\newblock {\em Comm. Anal. Geom.}, 3(3-4):597--644, 1995.

\bibitem[Gor97]{gordon:cmids}
Cameron~McA. Gordon.
\newblock Combinatorial methods in {D}ehn surgery.
\newblock In {\em Lectures at KNOTS '96 (Tokyo)}, volume~15 of {\em Ser. Knots
  Everything}, pages 263--290. World Sci. Publishing, River Edge, NJ, 1997.

\bibitem[Gor99]{gordon:ssadf}
Cameron~McA. Gordon.
\newblock Small surfaces and {D}ehn filling.
\newblock In {\em Proceedings of the Kirbyfest (Berkeley, CA, 1998)}, volume~2
  of {\em Geom. Topol. Monogr.}, pages 177--199 (electronic). Geom. Topol.
  Publ., Coventry, 1999.

\bibitem[Gri]{grisgsby:cdokfhocbc}
J.~Elisenda Grigsby.
\newblock {Combinatorial Description of Knot Floer Homology of Cyclic Branched
  Covers}, arXiv:math.GT/0610238.

\bibitem[GS99]{gompfstipsicz:4makc}
R.~Gompf and A.~Stipsicz.
\newblock {\em $4$-manifolds and Kirby calculus}, volume~20 of {\em Graduate
  Studies in Mathematics}.
\newblock American Mathematical Society, Providence, RI, 1999.

\bibitem[GT00]{gt:dsokwylsagok}
Hiroshi Goda and Masakazu Teragaito.
\newblock Dehn surgeries on knots which yield lens spaces and genera of knots.
\newblock {\em Math. Proc. Cambridge Philos. Soc.}, 129(3):501--515, 2000.

\bibitem[GW99]{gordonwu:taadf}
Cameron~McA. Gordon and Ying-Qing Wu.
\newblock Toroidal and annular {D}ehn fillings.
\newblock {\em Proc. London Math. Soc. (3)}, 78(3):662--700, 1999.

\bibitem[HMW92]{hmw:sotwlygsm}
Craig~D. Hodgson, G.~Robert Meyerhoff, and Jeffrey~R. Weeks.
\newblock Surgeries on the {W}hitehead link yield geometrically similar
  manifolds.
\newblock In {\em Topology '90 (Columbus, OH, 1990)}, volume~1 of {\em Ohio
  State Univ. Math. Res. Inst. Publ.}, pages 195--206. de Gruyter, Berlin,
  1992.

\bibitem[Lic62]{lickorish}
W.~B.~R. Lickorish.
\newblock A representation of orientable combinatorial {$3$}-manifolds.
\newblock {\em Ann. of Math. (2)}, 76:531--540, 1962.

\bibitem[Mor89]{morimoto:gofkils}
Kanji Morimoto.
\newblock Genus one fibered knots in lens spaces.
\newblock {\em J. Math. Soc. Japan}, 41(1):81--96, 1989.

\bibitem[MOS]{mos:acdokfh}
Ciprian Manolescu, Peter Ozsvath, and Sucharit Sarkar.
\newblock {A combinatorial description of knot Floer homology},
  arXiv:math.GT/0607691.

\bibitem[MOST]{most:oclfh}
Ciprian Manolescu, Peter Ozsvath, Zoltan Szabo, and Dylan Thurston.
\newblock {On combinatorial link Floer homology}, arXiv:math.GT/0610559.

\bibitem[MP]{mp:dfotm3m}
Bruno Martelli and Carlo Petronio.
\newblock {Dehn filling of the "magic" 3-manifold}, arXiv:math.GT/0204228.

\bibitem[OS]{os:kfhars}
Peter Ozsv{\'a}th and Zolt{\'a}n Szab{\'o}.
\newblock {Knot Floer homology and rational surgeries}, arXiv:math.GT/0504404.

\bibitem[OS05]{os:kfh+lss}
Peter Ozsv{\'a}th and Zolt{\'a}n Szab{\'o}.
\newblock On knot {F}loer homology and lens space surgeries.
\newblock {\em Topology}, 44(6):1281--1300, 2005.

\bibitem[Rie00]{rieck:hsomitdfs}
Yo'av Rieck.
\newblock Heegaard structures of manifolds in the {D}ehn filling space.
\newblock {\em Topology}, 39(3):619--641, 2000.

\bibitem[Sta62]{stallings:ofc3m}
John Stallings.
\newblock On fibering certain {$3$}-manifolds.
\newblock In {\em Topology of 3-manifolds and related topics (Proc. The Univ.
  of Georgia Institute, 1961)}, pages 95--100. Prentice-Hall, Englewood Cliffs,
  N.J., 1962.

\end{thebibliography}
\end{document}